\documentclass{elsarticle}
\usepackage{amsmath}
\usepackage{mathdots}
\usepackage{hyperref}
\usepackage{latexsym,amsfonts,amssymb}
\usepackage{bbm}
\usepackage[top=3.2cm,bottom=2.2cm,left=2.5cm,right=2.5cm]{geometry}
\usepackage{amsthm}
\usepackage{color}
\usepackage{hyperref}
\usepackage{comment}

\usepackage{soul} 
\soulregister\ref7 
\soulregister\cite7 
\usepackage{enumitem}

\hypersetup{
  colorlinks=true,
  citecolor=blue,
  linkcolor=red,
  urlcolor=green}

\allowdisplaybreaks[4]

\newtheorem{thm}{\quad \ \ Theorem}[section]
\newtheorem{cor}[thm]{\quad \ \ Corollary}

\newtheorem{lem}[thm]{\quad \ \ Lemma}

\numberwithin{equation}{section}

\begin{document}
\begin{frontmatter}

\title{Weighted moments and sharp convergence rates for supercritical multitype branching processes}
\author{Sheng Xiao$^{a,c}$\corref{myauthor}}
\cortext[mycorrespondingauthor]{Corresponding author}
\ead{xiaos\_2019@163.com}
\author{Yingqiu Li$^{a,c}$}
\ead{liyq-2001@163.com}
\author{Xingang Liang$^{b\ast}$}
\ead{liangxingang@th.btbu.edu.cn}

\address{$^{a}$School of Mathematics and Statistics, Hunan First Normal University, Changsha 410205, Hunan, China\\
$^{b}$School of Science, Beijing Technology and Business University, Beijing 100048, China\\
$^{c}$Hunan Provincial University Key Laboratory for Big Data Analysis and Application, Changsha 410205, Hunan, China
}

\begin{abstract}
Let $Z_n^i=(Z_n^i(1),\cdots,Z_n^i(d)), n\geqslant0$, be a supercritical $d$-type branching process starting from a single  particle of type $i\in\{1,\cdots,d\}$, where $d$ is a positive integer. The Kesten-Stigum theorem states that if the mean matrix $M$ of the offspring distribution has spectral radius $\rho>1$ with the associated positive left and right eigenvectors $u$ and $v$, then for all $1\leqslant i \leqslant d$, the sequence $W_n^i=\langle Z_n^i,u\rangle/\rho^n u(i)$ forms a non-negative martingale, and hence converges almost surely to a random variable $W^i$ as $n\rightarrow\infty$. In this paper,
we investigate the asymptotic properties of the martingale $\{W_n^i\}$, focusing on weighted moments and sharp rates of convergence. We establish necessary and sufficient conditions for the finiteness of $\phi$-moments of $W^i$, where $\phi(x)=x^{\alpha}\ell(x)$ with  $\alpha\geqslant1$ and $\ell$ is a positive function slowly varying at infinity.
We also determine the sharp exponential rate of convergence of $W_n^i$ to $W^i$ in the $\phi$-moment sense, identify a phase transition at the critical rate $\rho^{1/2},$ and characterize the behavior at the critical value. Our proofs rely on the double martingale structure of the fundamental martingale and on convex inequalities for martingales.
\end{abstract}

\end{frontmatter}

\section{Introduction and main results}
We consider a discrete time, discrete state multitype branching process (MBP) with a positive integer
$d$ of particle types. This model is a fundamental generalization of the classical Galton-Watson process and has attracted sustained interest due to its wide ranging applications in biology, physics, and epidemiology. Foundational results can be found in \cite{Harris1963,Mode1971,AN1972}. In recent decades, it has attracted considerable attention from many authors; see, for example, \cite{Hoppe1976,AV1995,MT2005,YY2009,ASY2023,TZ2025}. We now continue this line of investigation and provide some asymptotic properties of the martingale induced by this model.

As usual, we set $\mathbb{N}=\{0,1, \cdots\}$ be the set of non-negative integers, and $\mathbb{R}^d$ be the $d$-dimensional space of vectors with real coordinates, equipped with the scalar product and the $L^1$-norm defined respectively by
$$
\langle x, y\rangle:=\sum_{i=1}^d x(i) y(i) \quad \text { and } \quad\|x\|:=\sum_{i=1}^d|x(i)|, \quad x, y \in \mathbb{R}^d.
$$
We associate $d$ probability generating functions: for each $1 \leqslant r \leqslant d$,
$$
f^r(s)=\sum_{k_1, \cdots, k_d=0}^{\infty} p_{k_1, \cdots, k_d}^r s_1^{k_1} \cdots s_d^{k_d}, \quad s=\left(s_1, \cdots, s_d\right) \in[0,1]^d,
$$
where $p_{k_1, \cdots, k_d}^r$ denotes the probability that a particle of type $r$ produces $k_1$ offspring of type 1,$\cdots$, and $k_d$ offspring of type $d$. A MBP $\left(Z_n\right)$ is a process with values in $\mathbb{N}^d$ such that, for all $n \geqslant 0$,
\begin{align}\label{eq1}
Z_{n+1}=\sum_{r=1}^d \sum_{l=1}^{Z_n(r)} N_{l, n}^r,
\end{align}
where $Z_0 \in \mathbb{N}^d$ is fixed, $Z_n(j)$ represents the number of particles of type $j$ in generation $n$, and $N_{l, n}^r(j)$ denotes the number of offspring of type $j$ at time $n+1$ produced by the $l$-th particle of type $r$ in generation $n$. The random vectors $N_{l, n}^r= (N_{l, n}^r(1), \cdots, N_{l, n}^r(d))$, indexed by $l \geqslant 1, n \geqslant 0$, $1\leqslant r \leqslant d$, are independent and have the same probability generating function $f^r$. When the process starts from a single initial particle of type $i$, that is, $Z_0=e_i$, we write $Z_n^i$ instead of $Z_n$. Here $e_i$ is the $i$-th standard unit vector, where the $i$-th component is $1$ and all other components are $0$.

Denote the mean matrix of the offspring distributions by $M$, where the $(i,j)$-th entry $M(i,j)$ of $M$ represents the expected number of type $j$ offspring produced in the first generation by a single particle of type $i$. Let $M^n$ denote the $n$-th power of the matrix $M$, and let $M^n(i,j)$ denote its $(i,j)$-th entry. We always assume that there exists an $n$ such that $M^n(i,j)>0$ for all $1\leqslant i, j\leqslant d$; in this case, the process $(Z_n)$ is said to be positive regular.
According to the Perron-Frobenius theorem, the mean matrix $M$ has the maximal eigenvalue $\rho>0$, and the associated positive left and right eigenvectors $u=(u(1), \cdots, u(d))$ and $v=(v(1), \cdots, v(d))$ with the normalization $\|u\|=1$ and $\langle v, u\rangle=1$. That means
\begin{align}\label{eq4}
Mu=\rho u;
\end{align}
and moreover, by induction,
\begin{align}\label{eq3}
M^nu=\rho^nu.
\end{align}

In this paper, we always assume that $\rho>1$, which means that the branching process is in the supercritical regime. The classical Kesten-Stigum theorem (see \cite{KS1966}) shows that for all $1 \leqslant i, j \leqslant d$, as $n \rightarrow+\infty$,
\begin{align}\label{eq2}
\frac{Z_n^i(j)}{\rho^n v(j)} \rightarrow W^i u(i) \quad \text {a.s.,}
\end{align}
where $W^i$ is a non-negative random variable,
 and is nondegenerate for all $i$ if and only if
$$\mathbb{E}\left(Z_1^i(j) \log ^{+} Z_1^i(j)\right)<+\infty \ \ \forall1\leqslant i,j\leqslant d;$$ moreover, when it is nondegenerate, $\mathbb{E} W^i=1$. The proof of (\ref{eq2}) is based on the fundamental non-negative martingale
\begin{align}\label{eq6}
W_0^i=1, \quad W_n^i=\frac{\left\langle Z_n^i, u\right\rangle}{\rho^n u(i)}, \quad n \geqslant 1,
\end{align}
which converges a.s. to $W^i$. When $W^i$ is nondegenerate, the study of its moments becomes a natural and important problem. The case $d=1$ has been extensively investigated; see, for example, \cite{Harris1963,AN1972,BD1974,AR2004,LL2013}. We briefly recall that, in this setting, $(W_n^i)$ reduces to the fundamental martingale $(Z_n/m^n)$ of a Galton-Watson process $(Z_n)$, where
$m$ denotes the mean of the offspring distribution. For the limit $W$
of the sequence $(Z_n/m^n)$, Bingham and Doney \cite{BD1974} proved the equivalence
\begin{align}\label{eq23}
\mathbb{E}(Z_1)^\alpha<\infty\Longleftrightarrow\mathbb{E}W^\alpha<\infty \ \  \text{for} \ \  \alpha>1.
\end{align}
A natural and interesting question is whether this equivalence remains valid under functional transformations. Using a deep Tauberian theorem, Bingham and Doney \cite{BD1974} showed that when $\alpha>1$ is not an integer, $\mathbb{E} W^\alpha \ell(W)<\infty$ if and only if $\mathbb{E}Z_1^\alpha \ell\left(Z_1\right)<\infty$, where $\ell$ is a positive function slowly varying at infinity. By a nice martingale argument, Alsmeyer and R\"{o}sler \cite{AR2004} showed that the equivalence remains true when $\alpha$ is not of the form $2^n$ for some integer $n \geqslant 1$. Liang and Liu \cite{LL2013} further refined this martingale approach and proved that the equivalence holds for all $\alpha>1$. More precisely, they showed that the following assertions are equivalent:
\begin{align}\label{eq24}
(i) \ \mathbb{E}(Z_1)^\alpha\ell(Z_1)<\infty;\ \ (ii) \ \mathbb{E}(\widehat{W})^\alpha\ell(\widehat{W})<\infty; \ \ (iii) \ \mathbb{E}(W)^\alpha\ell(W)<\infty,
\end{align}
where $\widehat{W}=\sup_{n\geqslant0}W_n$. This
is the so-called comparison theorem for the weighted moments of $Z_1$, $\widehat{W}$, and $W$.  In fact, the study of weighted moments in branching systems has attracted sustained attention, see for instance, \cite{LHL2011}, \cite{LL2013a} and \cite{WLLL2019a} for investigations in branching processes in random environments, and \cite{AK2010}, \cite{LL2015}, \cite{WH2019}, \cite{WLLL2019b} for studies on branching random walks.

For a MBP with $d>1$, there are few results on the existence of weighted moments of the limit $W^i$. One of the aims of this paper is therefore to investigate this question and to establish criteria for the finiteness of weighted moments of the form $\mathbb{E}(W^i)^\alpha \ell(W^i),$ where $\alpha\geqslant1$ and $\ell$ is slowly varying at infinity. To this end, we extend the martingale methods developed in \cite{AR2004} and \cite{LL2013} to the multitype setting. Beyond the existence of weighted moments, a natural and more refined question is to determine how fast the fundamental martingale $W_n^i$ converges to its limit $W^i$ in the weighted moment sense. More precisely, for $\phi(x)=x^\alpha\ell(x)$, we are interested in the asymptotic behavior of
$\mathbb{E}\phi(\delta^n|W^i-W_n^i|)$
as $n\to\infty$,
together with the associated phase transition and the behavior at the critical value.
To the best of our knowledge, such weighted moment convergence rates have not been systematically investigated, even for the classical Galton-Watson process.
A key tool in our analysis is the double martingale structure underlying $(W_n^i)$, which, together with suitable convex inequalities for martingales, provides an effective framework for studying weighted behavior.

Our first result concerns the uniform boundedness of the $\alpha$-th moments of $(W_n^i)$, which is an extension of (\ref{eq23}) to the case $d>1$. This result is instrumental in the proof of the existence of weighted moments of $W^i$.
\begin{thm}\label{Thm1}
Let $\alpha>1$. The following assertions are equivalent:
\begin{itemize}[leftmargin=3em]
    \item [(a)] $\mathbb{E}\left(Z^i_1(j)\right)^{\alpha}<\infty \ \ \forall1\leqslant i,j\leqslant d;$
    \item [(b)] $\sup_{n\geqslant0}\mathbb{E}(W_{n}^i)^{\alpha}<\infty \ \ \forall1\leqslant i\leqslant d;$
    \item[(c)] $W_n^i\longrightarrow W^i \ in \ L^\alpha \ \forall1\leqslant i\leqslant d;$
    \item[(d)]  $0<\mathbb{E}(W^i)^\alpha<\infty \ \ \forall1\leqslant i\leqslant d.$
\end{itemize}
\end{thm}

Under the Furstenberg-Kesten condition (see \cite{FK1960}), Grama, Liu and Pin \cite{GLP2022} established necessary and sufficient conditions for the $L^p$ convergence of the fundamental martingale associated with a MBP in random environments. When the environmental sequence is constant, their result reduces to the corresponding statement for a MBP in deterministic environments. Their approach relies on an ingenious environmental shift technique to estimate the quenched $\alpha$-th moment of $W_{n+1}^i-W_n^i$. Our method is quite different:
we exploit the double martingale structure of $(W_n^i)$ together with Burkholder's inequalities to derive bounds for $\mathbb{E}|W_{n+1}^i- W_n^i|^\alpha.$

Theorem \ref{Thm1} indicates that $Z^i_1(j)$ and $W^i$ have similar tail behavior. In the following, we will establish comparison theorems between weighted moments of $Z^i_1(j)$ and $W^i,$ thereby extending (\ref{eq24}) to the case $d>1.$

Recall that a positive and measurable function $\ell$ defined on $[0, \infty)$ is called slowly varying at infinity if $\lim _{x \rightarrow \infty} \ell(t x) / \ell(x)=1$ for all $t>0$. By the representation theorem (see \cite{BGT1987}), any slowly varying function $\ell$ is of the form
$$
\ell(x)=c(x) \exp \left(\int_{a_0}^x \epsilon(t) \frac{d t}{t}\right), \quad x>a_0,
$$
where $a_0 \geqslant 0, c(\cdot)$ and $\epsilon(\cdot)$ are measurable with $c(x) \rightarrow c$ for some constant $c \in(0, \infty)$ and $\epsilon(x) \rightarrow 0$ as $x \rightarrow \infty$. The value of $a_0$ and those of $\ell(x)$ on $\left[0, a_0\right]$ will not be important; we always assume that $\ell$ is bounded on compact sets of $[0, \infty)$. For convenience, we often take $a_0=1$. Moreover, the function $c(x)$ in the representation of $\ell(x)$ has no influence on the finiteness of the moments, so that we can suppose without loss of generality that $c(x)=1.$
Set
$$\widehat{W}^i:=\sup_{n\geqslant0}W_n^i, \ \forall 1\leqslant i\leqslant d.$$

When $\alpha>1$, we have the following comparison theorem between weighted moments of $Z^i_1(j)$ and $W^i$:
\begin{thm}\label{Thm2}
Let $\alpha>1$ and $\ell:[0, \infty) \mapsto[0, \infty)$ be a function slowly varying at $\infty$. Let $\phi(x)=x^\alpha\ell(x)$. The following assertions are equivalent:
\begin{itemize}[leftmargin=3em]
  \item [(a)] $\mathbb{E}\phi\left(Z^i_1(j)\right)<\infty \ \ \forall1\leqslant i,j\leqslant d;$
  \item [(b)] $\mathbb{E}\phi(W_1^{i})<\infty \ \forall1\leqslant i\leqslant d;$
  \item [(c)] $\mathbb{E}\phi(\widehat{W}^i)<\infty \ \ \forall1\leqslant i\leqslant d;$
  \item [(d)] $\mathbb{E}W^i=1$ \ and \ $\mathbb{E}\phi(W^i)<\infty, \ \ \forall1\leqslant i\leqslant d.$
\end{itemize}
\end{thm}
We also consider the case $\alpha=1$, where the situation is different. For a measurable function $\ell:[0,\infty) \mapsto [0,\infty)$, set
\begin{equation*}
    \tilde{\ell}(x)=
    \begin{cases}
        \int_1^x \ell(t)\frac{dt}{t}, &\quad \text{if}\ x>1,\\
        0, &\quad \text{if}\ x\leqslant 1.
    \end{cases}
\end{equation*}
Then when $\alpha=1$ and $\ell$ is concave, we have the following comparison theorem:
\begin{thm}\label{Thm3}
    Let $\ell: [0,\infty) \mapsto [0,\infty)$ be a positive increasing function which is concave on $[0,\infty)$.
    If $\mathbb{E} W_1^i \tilde{\ell}(W_1^i)<\infty$ for $1\leqslant i \leqslant d$, then
    \begin{equation*}
        \mathbb{E} \widehat{W}^i \ell(\widehat{W}^i)<\infty \quad \text{and}\quad \mathbb{E} W^i \ell(W^i)<\infty, \ \ \forall \ 1\leqslant i\leqslant d.
    \end{equation*}
\end{thm}
As a corollary of Theorem \ref{Thm3}, we have:
\begin{cor} \label{Cor1}
    Let $\ell:[0,\infty)\mapsto[0,\infty)$ be non-decreasing and slowly varying at $\infty$, such that $\ell(x)=\tilde{\ell}_0(x)$ for some slowly varying function $\ell_0(x)$.
    If $\mathbb{E} W_1^i \tilde{\ell}(W_1^i)<\infty$ for $1\leqslant i \leqslant d$, then
    \begin{equation*}
        \mathbb{E} \widehat{W}^i \ell(\widehat{W}^i)<\infty \quad \text{and}\quad \mathbb{E} W^i \ell(W^i)<\infty, \ \ \forall \ 1\leqslant i\leqslant d.
    \end{equation*}
\end{cor}

Corollary \ref{Cor1} is a direct consequence of Theorem \ref{Thm3}. Since $\ell'(x)= x^{-1}\ell_0(x)$ is a regularly varying function of order $-1$, we have that $\ell'(x)$ is equivalent to $\psi(x) = \inf\{\ell'(t):1\leqslant t\leqslant x\} $ when $x\rightarrow\infty$, where $\psi(x)$ is non-increasing; moreover, without loss of generality, we take $\ell(x)=\int_1^x\psi(t)dt$ ($x\geqslant 1$), which is a concave function. Therefore, the assumptions in Theorem \ref{Thm3} hold true under the assumptions on $\ell(x)$ in Corollary \ref{Cor1}, and hence we obtain the same conclusion as that in Theorem \ref{Thm3}.

We now turn from the existence of weighted moments to the rate of convergence of $W_n^i$ to $W^i$ in the weighted moment sense. More precisely, we investigate the asymptotic behavior of
$\mathbb E\phi(\delta^n|W^i-W_n^i|).$ The following theorem provides upper estimates for this quantity under different moment assumptions and identifies the admissible exponential rates.

In the following results, we will always assume that $\mathbb{P}(W_1^j=1)<1$ for some $1\leqslant j\leqslant d,$ which excludes the trivial case where $W^i=1$ a.s. for all $1\leqslant i\leqslant d.$

\begin{thm}\label{Thm6}
Let $\alpha>1$ and $\ell:[0, \infty) \mapsto[0, \infty)$ be a function slowly varying at $\infty$. Let $\phi(x)=x^\alpha\ell(x)$.
\begin{itemize}[leftmargin=3em]
  \item[(i)] Suppose that
$$\mathbb{E}\phi(Z_1^i(j))<\infty \ \  \forall 1\leqslant i,j\leqslant d.$$
Then for every $1<\delta<\rho^{1-1/\min{\{\alpha,2\}}}$, we have
\begin{equation}\label{10eq1}
\lim_{n\to\infty}\mathbb{E}\phi(\delta^n|W^i - W_n^i|)=0 \ \  \forall 1\leqslant i\leqslant d.
\end{equation}
\item[(ii)] Let $\alpha\in(1,2)$. Suppose that there exists some $\eta\in(\alpha,2]$ such that
$$\mathbb{E}(Z_1^i(j))^{\eta}<\infty \ \  \forall 1\leqslant i,j\leqslant d.$$
Then \eqref{10eq1} holds for every $1<\delta<\rho^{1-1/\eta}.$
\end{itemize}
\end{thm}

We next turn to the behavior at the critical rate $\rho^{1/2}$. When $\alpha\neq2$, the second moment condition, together with the weighted moment assumption, is sufficient to obtain the uniform upper bound at the critical value. The case $\alpha=2$ is more delicate, since it lies at the boundary between the convexity and concavity regimes arising in the martingale inequalities. Nevertheless, the same critical upper bound remains valid under some additional regularity conditions on the weighted moment function.

\begin{thm}\label{Thm6-1}
Let $\alpha>1$ and $\ell:[0, \infty) \mapsto[0, \infty)$ be a function slowly varying at $\infty$. Let $\phi(x)=x^\alpha\ell(x)$.
Suppose that
$$\mathbb{E}(Z_1^i(j))^{2}<\infty \ \ \text{and} \ \ \mathbb{E}\phi(Z_1^i(j))<\infty \ \  \forall 1\leqslant i,j\leqslant d.$$
\begin{itemize}[leftmargin=3em]
  \item[(i)] If $\alpha\neq2$, then, at the critical value $\delta=\rho^{1/2}$, we have
\begin{equation}\label{10eq6}
\limsup_{n\to\infty}\mathbb{E}\phi(\delta^n|W^i - W_n^i|)<\infty \ \  \forall 1\leqslant i\leqslant d.
\end{equation}
  \item[(ii)] For $\alpha=2$, if one of the following conditions is satisfied:

  (a) the function $x\rightarrow \phi(x^{1/2})$ is either convex or concave on $[0, \infty);$

  (b) there exists a constant $C>0$ such that $\ell(tx)\leqslant C\ell(x), \ \forall 0<t\leqslant1,x>0,$

  then \eqref{10eq6} holds.
\end{itemize}
\end{thm}

Theorem \ref{Thm6-1} provides upper estimates at the critical rate. We next establish the corresponding lower estimates. The following result shows that $\rho^{1/2}$ is indeed the critical exponential scale: above this threshold, the weighted moments diverge to infinity, whereas at the threshold they remain bounded away from zero.

\begin{thm}\label{Thm7}
Let $\alpha>1$ and $\ell:[0, \infty) \mapsto[0, \infty)$ be a function slowly varying at $\infty$. Let $\phi(x)=x^\alpha\ell(x)$.  
Then, for all $1 \leqslant i \leqslant d$, the following assertions hold:
\begin{itemize}[leftmargin=3em]
\item[(i)] If $\delta>\rho^{1/2},$ then
\begin{equation*}
\lim_{n\to\infty}\mathbb{E}\phi(\delta^n|W^i - W_n^i|) =\infty.
\end{equation*}
\item[(ii)] At the critical value $\delta=\rho^{1/2},$  we have
\begin{equation*}
\liminf_{n\to\infty}\mathbb{E}\phi(\delta^n|W^i - W_n^i|)>0.
\end{equation*}
\end{itemize}
\end{thm}

Combining Theorems \ref{Thm6} and \ref{Thm7}, we obtain a sharp phase transition result under suitable moment conditions. In particular, the critical exponential rate is $\rho^{1/2}$: below this rate the weighted moments converge to zero, while above it they diverge to infinity.

\begin{cor}
Let $\alpha>1$ and $\ell:[0, \infty) \mapsto[0, \infty)$ be a function slowly varying at $\infty$. Let $\phi(x)=x^\alpha\ell(x)$.
Suppose that
$$\mathbb{E}(Z_1^i(j))^{2}<\infty \ \ \text{and} \ \ \mathbb{E}\phi(Z_1^i(j))<\infty \ \  \forall 1\leqslant i,j\leqslant d.$$
Then, for all $1 \leqslant i \leqslant d$,
\begin{equation*}
\lim_{n\to\infty}\mathbb{E}\phi(\delta^n|W^i - W_n^i|) =
\begin{cases}
0, & 1 < \delta < \rho^{1/2}, \\
\infty, & \delta > \rho^{1/2}.
\end{cases}
\end{equation*}
In particular, $\mathbb{E}\phi\left(\delta^n|W^i - W_n^i|\right) \longrightarrow 0$
if and only if $1 < \delta < \rho^{1/2}.$
\end{cor}

The preceding corollary contains the corresponding $L^\alpha$-convergence result as a special case. In fact, taking $\ell$ to be any positive constant, we have $\phi(x)=cx^\alpha$ for some $c>0$. Consequently, the weighted moment phase transition yields the following sharp exponential $L^\alpha$-convergence rate, including the behavior at the critical value.
\begin{cor}
Let $\alpha>1$. Suppose that
\begin{align}\label{Xiaoeq10}
\mathbb{E}(Z_1^i(j))^{\alpha\vee2}<\infty \ \  \forall 1\leqslant i,j\leqslant d.
\end{align}
Then, for all $1 \leqslant i \leqslant d$,
\begin{equation*}
\lim_{n\to\infty} \delta^n \left(\mathbb{E}|W^i - W_n^i|^\alpha\right)^{1/\alpha} =
\begin{cases}
0, & 1 < \delta < \rho^{1/2}, \\
\infty, & \delta > \rho^{1/2}.
\end{cases}
\end{equation*}
In particular, $\delta^n \left(\mathbb{E}|W^i - W_n^i|^\alpha\right)^{1/\alpha} \longrightarrow 0$
if and only if $1 < \delta < \rho^{1/2}.$

Moreover, at the critical value $\delta = \rho^{1/2}$,
\begin{equation*}
0 < \liminf_{n\to\infty} \delta^n \left(\mathbb{E}|W^i - W_n^i|^\alpha\right)^{1/\alpha} \leqslant \limsup_{n\to\infty} \delta^n \left(\mathbb{E}|W^i - W_n^i|^\alpha\right)^{1/\alpha} < \infty.
\end{equation*}
\end{cor}

The remainder of this paper is organized as follows. In Section 2, we explain the double martingale structure of $(W_n^i)$ and present the proof of Theorem \ref{Thm1}. Sections 3 and 4 are devoted to the proofs of Theorems \ref{Thm2} and \ref{Thm3}, respectively, which concern weighted moments. Sections 5-7 establish the convergence rates by proving Theorems 1.5, 1.6, and 1.7, respectively.

Throughout the paper, the symbol $C$, possibly with some indices, always denotes a positive constant whose value may differ from line to line.

\section{Double martingale structure and Proof of Theorem \ref{Thm1}}
\subsection{Double martingale structure}

To obtain the desired results, we introduce a double martingale structure on $(W_n^i)$, which will be used repeatedly in the proofs of our main results. 
Let
$$\mathcal{F}_0=\{\varnothing,\Omega\}, \ \ \mathcal{F}_n=\sigma\{N_{l,k}^r:0\leqslant k\leqslant n-1, 1\leqslant r\leqslant d,l\geqslant1\}.$$
Clearly $(W_{n}^i)$ is adapted to $(\mathcal{F}_n)$. Noting that $N_{l,n}^r$ is independent of $\mathcal{F}_n$ and has the same distribution as $N_{1,0}^r$, it follows from the definitions of $Z_n^i$ and $M$ that, for all $n\geqslant0$ and $1\leqslant i\leqslant d$,
\begin{equation*}
\mathbb{E}(Z_{n+1}^i |\mathcal{F}_n)=\sum_{r=1}^d\sum_{l=1}^{Z_n^i(r)} \mathbb{E} N_{l,n}^r = \sum_{r=1}^d Z_n^i(r) \mathbb{E}N_{1,0}^r=M^T Z_n^i.
\end{equation*}
Moreover, using \eqref{eq4}, we obtain that for all $n\geqslant 0$ and $1\leqslant i\leqslant d$,
\begin{equation*}
\mathbb{E}(W_{n+1}^i|\mathcal{F}_n)=\frac{\langle \mathbb{E}(Z_{n+1}^i|\mathcal{F}_n), u \rangle}{\rho^{n+1} u(i)}=\frac{\langle M^T Z_n^i, u \rangle}{\rho^{n+1}u(i)}=\frac{\langle Z_n^i, Mu\rangle}{\rho^{n+1}u(i)}=\frac{\langle Z_n^i, \rho u \rangle}{\rho^{n+1}u(i)}=W_n^i,
\end{equation*}
which means that $(W_n^i)_{n\geqslant 0}$ forms a martingale with respect to the filtration $(\mathcal{F}_n)_{n\geqslant 0}$.

Fix $n\geqslant0$ and let
\begin{align}\label{eq5}
D_{n,0}^i=0, \ \ D_{n,k}^i=\sum_{r=1}^d \frac{u(r)}{\rho^n u(i)} \sum_{l=1}^{k}\left(W_{l, n}^r-1\right)\mathbf{1}_{\{Z_n^i(r)\geqslant l\}},
\end{align}
where
$$
W_{l, n}^r:=\frac{\langle N_{l, n}^r, u\rangle}{\rho u(r)}.
$$
Let $\mathcal{G}_{n,0}=\mathcal{F}_n, \ \mathcal{G}_{n,k}=\sigma\{\mathcal{F}_n, W_{l,n}^r: 1\leqslant l\leqslant k,1\leqslant r\leqslant d\}.$ It is not difficult to verify that the sequence $(D_{n,k}^i)_{k\geqslant0}$ forms a martingale with respect to the filtration $(\mathcal{G}_{n,k})_{k\geqslant0}$ .
In fact, as $N_{k,n}^r$ is independent of the filtration $\mathcal{G}_{n,k-1}$ with mean vector $M(r,\cdot)$, we have that
for all $k\geqslant1$ and $1 \leqslant i \leqslant d,$
\begin{align*}
\mathbb{E}\left(W_{k, n}^r|\mathcal{G}_{n,k-1}\right)
&=\mathbb{E}\left(\frac{\langle N_{k, n}^r, u\rangle}{\rho u(r)}|\mathcal{G}_{n,k-1}\right)
=\frac{\langle \mathbb{E}N_{k,n}^r,u\rangle}{\rho u(r)}
=\frac{\langle M(r,\cdot),u\rangle}{\rho u(r)}
=\frac{(Mu)(r)}{\rho u(r)} =1;
\end{align*}
moreover,
\begin{align*}
\mathbb{E}(D_{n,k}^i|\mathcal{G}_{n,k-1})&=\mathbb{E}\left(\sum_{r=1}^d \frac{u(r)}{\rho^n u(i)} \sum_{l=1}^{k}\left(W_{l, n}^r-1\right)\mathbf{1}_{\{Z_n^i(r)\geqslant l\}}|\mathcal{G}_{n,k-1}\right)\\
&=D_{n,k-1}^i+\sum_{r=1}^d\frac{u(r)}{\rho^n u(i)}\mathbf{1}_{\{Z_n^i(r)\geqslant k\}}\mathbb{E}\left(\left(W_{k, n}^r-1\right)|\mathcal{G}_{n,k-1}\right)\\
&=D_{n,k-1}^i.
\end{align*}

Recalling the definition of $Z_n^i$, we have
\begin{equation*}
\langle Z_{n+1}^i/\rho-Z_n^i,u \rangle
=\sum_{r=1}^d \sum_{l=1}^{Z_n^i(r)}\left(\left\langle\frac{N_{l,n}^r}{\rho},u \right\rangle -u(r)\right)
=\sum_{r=1}^d u(r)\sum_{l=1}^{Z_n^i(r)}(W_{l,n}^r-1).
\end{equation*}
Furthermore, it follows from the definition of $W_n^i$ that
\begin{eqnarray}\label{eq2-liang}
W_{n+1}^i-W_n^i
&=&\frac{\langle Z_{n+1}^i/\rho-Z_n^i,u\rangle}{\rho^n u(i)}
=\sum_{r=1}^d \frac{u(r)}{\rho^n u(i)}\sum_{l=1}^{Z_n^i(r)}(W_{l,n}^r-1)\notag\\
&=& \sum_{r=1}^d \frac{u(r)}{\rho^n u(i)}\sum_{l=1}^{\infty}(W_{l,n}^r-1)\mathbf{1}_{\{Z_n^i(r)\geqslant l\}}
= \lim_{k\rightarrow\infty}D_{n,k}^i.
\end{eqnarray}
Therefore, $(W_{n+1}^i-W_n^i)_{n\geqslant0}$ is a martingale difference sequence with respect to the filtration $( \mathcal{F}_n)_{n\geqslant0}$, and  can be regarded as the limit of the martingale $(D_{n,k}^i,\mathcal{G}_{n,k})_{k\geqslant0}$; that is why we refer to it as the double martingale structure.
Our proofs of  the theorems rely mainly  on this double martingale structure, and  we will  apply classical convex inequalities to the martingale difference sequences $(D_{n,k}^i-D_{n,k-1}^i)$ and $(W_{n+1}^i-W_n^i).$

\subsection{Proof of Theorem \ref{Thm1}}
Before giving the proof of Theorem \ref{Thm1}, we first recall the following Burkholder's inequality for martingales.
\begin{lem} (Burkholder's inequality, see e.g. \cite{CT2012}).
Let $\{f_n\}$ be a $L^1$ martingale with $f_0=0$. Let $Q_n=(\sum_{k=1}^n(f_k-f_{k-1})^2)^{1 / 2}$ and
$Q=(\sum_{k=1}^{\infty}(f_k-f_{k-1})^2)^{1 / 2}$.
Then for all $p>1,$ there exist two positive constants $c_p=(p-1) / 18 p^{3 / 2}, C_p=18 p^{3 / 2} /(p-1)^{1 / 2}$ such that
$$
\begin{gathered}
c_p\left\|Q_n\right\|_p \leqslant\left\|f_n\right\|_p \leqslant C_p\left\|Q_n\right\|_p, \\
c_p\|Q\|_p \leqslant \sup _n\left\|f_n\right\|_p \leqslant C_p\|Q\|_p.
\end{gathered}
$$
\end{lem}
The following
lemma is a direct consequence of the Marcinkiewicz-Zygmund inequality.
\begin{lem}\label{lem4}(\cite{Liu2001}, Lemma 1.4).
Let $(X_k)_{k\geqslant 1}$ be a sequence of i.i.d. random centered variables. Then we have for $\alpha>1,$
$$
\mathbb{E}\left|\sum_{k=1}^n X_k\right|^\alpha \leqslant\left\{\begin{array}{ll}
\left(B_\alpha\right)^\alpha \mathbb{E}\left|X_k\right|^\alpha n, \quad &\text { if } 1<\alpha \leqslant 2, \\
\left(B_\alpha\right)^\alpha \mathbb{E}\left|X_k\right|^\alpha n^{\frac{\alpha}{2}}, \quad &\text { if } \alpha>2,
\end{array}\right.
$$
where $B_\alpha=2 \min \left\{k^{1 / 2}: k \in \mathbb{N}, k \geqslant \frac{\alpha}{2}\right\}$ is a constant depending only on $\alpha.$
\end{lem}
\begin{proof}[Proof of Theorem \ref{Thm1}.]
We first prove the implication $(a)\Rightarrow(b).$ Noting that the function $x\mapsto x^{\alpha}$ is convex (as $\alpha>1$) and the eigenvector $u$ is normalized with all components positive, we have by Jensen's inequality that for all $1\leqslant i\leqslant d$,
\begin{align*}
\mathbb{E}(W_{1}^i)^{\alpha}
&=\mathbb{E}\left(\sum_{j=1}^d\frac{Z_1^i(j)u(j)}{\rho u(i)}\right)^\alpha
\leqslant\mathbb{E}\sum_{j=1}^du(j)\left(\frac{Z_1^i(j)}{\rho u(i)}\right)^\alpha
=\sum_{j=1}^d\frac{u(j)}{(\rho u(i))^\alpha}\mathbb{E}\left(Z_1^i(j)\right)^\alpha<\infty.
\end{align*}
%
We divide the proof into two cases, where $\alpha\in (1,2]$ or $\alpha>2$.

For the case $\alpha\in(1,2]$, using Burkholder's inequality to $(W_{n+1}^i-W_{n}^i)$ and the subadditivity of the function $x\rightarrow x^{\alpha/2}$, we have that for all $1\leqslant i\leqslant d,$
\begin{align}\label{eq1-liang}
\sup_{n\geqslant0}\mathbb{E}|W^i_n-1|^\alpha&=\sup_{n\geqslant0}\mathbb{E}\left|\sum_{k=0}^n(W_{k+1}^i-W_k^i)\right|^\alpha
\leqslant C\mathbb{E}\left(\sum_{n=0}^{\infty}|W^i_{n+1}-W^i_{n}|^2\right)^{\alpha/2}\notag\\
&\leqslant C\mathbb{E}\sum_{n=0}^{\infty}|W^i_{n+1}-W^i_{n}|^\alpha
=C\sum_{n=0}^{\infty}\mathbb{E}|W^i_{n+1}-W^i_{n}|^\alpha.
\end{align}
Recall that under $\mathcal{F}_n$, $W_{l,n}^r-1$ are i.i.d. random centered variables and are independent of $Z_n^i(r)$, so by Lemma \ref{lem4}, we see that
\begin{eqnarray*}
\mathbb{E}\left(\left|\sum_{l=1}^{Z_n^i(r)} (W_{l,n}^r-1)\right|^{\alpha}|\mathcal{F}_n\right)
\leqslant B_{\alpha}^{\alpha}Z_n^i(r) \mathbb{E}|W_{l,n}^r-1|^{\alpha}
= B_{\alpha}^{\alpha}Z_n^i(r) \mathbb{E}|W_{1,0}^r-1|^{\alpha}.
\end{eqnarray*}
Moreover, applying the convexity of the function $x\mapsto x^{\alpha}$ (and $\sum_{j=1}^{d}u(j)=1$),
we get that for all $1\leqslant i\leqslant d,$
\begin{align}\label{eq7}
\mathbb{E}\left(\left|W_{n+1}^i-W_n^i\right|^{\alpha}|\mathcal{F}_n\right)
&=\mathbb{E}\left[\left|\sum_{r=1}^d \frac{u(r)}{\rho^n u(i)} \sum_{l=1}^{Z_n^i(r)}\left(W_{l, n}^r-1\right)\right|^{\alpha}|\mathcal{F}_n\right]\nonumber\\
&\leqslant \frac{1}{\left(\rho^n u(i)\right)^{\alpha}}\sum_{r=1}^d u(r)\mathbb{E}\left(\left|\sum_{l=1}^{Z_n^i(r)}\left(W_{l, n}^r-1\right)\right|^{\alpha}|\mathcal{F}_n\right)\nonumber\\
&\leqslant \frac{1}{\left(\rho^n u(i)\right)^{\alpha}}\sum_{r=1}^d u(r) B_{\alpha}^{\alpha}Z_{n}^{i}(r)\mathbb{E}\left|W_{1, 0}^r-1\right|^{\alpha}\nonumber\\
&\leqslant C \rho^{(1-\alpha)n}W_n^i,
\end{align}
where $C=B_{\alpha}^{\alpha}\left(u(i)\right)^{1-\alpha}\cdot \max\limits_{1\leqslant j\leqslant d} \mathbb{E}|W_{1,0}^j-1|^{\alpha}>0$ is finite, as $u(i)>0$ for all $1\leqslant i \leqslant d$ and
\begin{equation}\label{eq3-liang}
\max_{1\leqslant j\leqslant d} \mathbb{E}|W_{1,0}^j-1|^{\alpha}
=\max_{1\leqslant j\leqslant d} \mathbb{E}|W_{1}^j-1|^{\alpha}
\leqslant \max_{1\leqslant j\leqslant d} 2\left(\mathbb{E}(W_{1}^j)^{\alpha}+1\right)
<\infty\,.
\end{equation}
%
%
%
Then taking expectation on both sides in \eqref{eq7}, as $(W_n^i)_{n\geqslant 0}$ is a martingale with $W_0^i=1$, we obtain
\begin{align*}
\mathbb{E}\left|W_{n+1}^i-W_n^i\right|^{\alpha}
\leqslant C\rho^{(1-\alpha)n}\mathbb{E} W_n^i =C\rho^{(1-\alpha)n}.
\end{align*}
%
Therefore,  it follows from the supercritical assumption $\rho>1$ and \eqref{eq1-liang}
that for all $1\leqslant i\leqslant d,$
\begin{align}\label{eq6-liang}
\sup_{n\geqslant0}\mathbb{E}|W^i_n-1|^\alpha&\leqslant C\sum_{n=0}^\infty \rho^{(1-\alpha)n}<\infty\,.
\end{align}

We now turn to the case $\alpha>2$. Using Burkholder's inequality to the martingale difference
sequence $(W^i_{n+1}-W^i_{n})$, we have that for all $1\leqslant i\leqslant d,$
\begin{align*}
\sup_{n\geqslant0}\mathbb{E}|W^i_n-1|^\alpha
\leqslant C\mathbb{E}\left(\sum_{n=0}^{\infty}|W^i_{n+1}-W^i_{n}|^2\right)^{\alpha/2};
\end{align*}
moreover, by the triangular inequality in $L^{\alpha/2},$
\begin{align*}
\left[\mathbb{E}\left(\sum_{n=0}^{\infty}|W^i_{n+1}-W^i_{n}|^2\right)^{\alpha/2}\right]^{2/\alpha}
\leqslant\sum_{n=0}^{\infty}\left(\mathbb{E}|W^i_{n+1}-W^i_{n}|^{\alpha}\right)^{2/\alpha}.
\end{align*}
Then
\begin{align}\label{eq10}
\sup_{n\geqslant0}\mathbb{E}|W^i_n-1|^\alpha\leqslant C\left[\sum_{n=0}^{\infty}\left(\mathbb{E}|W^i_{n+1}-W^i_{n}|^{\alpha}\right)^{2/\alpha}\right]^{\alpha/2}.
\end{align}

In the following, we will deal with the term $\mathbb{E}|W^i_{n+1}-W^i_{n}|^{\alpha}$. Recalling that $W_{n+1}^i-W_n^i=\lim\limits_{k\rightarrow\infty}D_{n,k}^i$ (cf. \eqref{eq2-liang}),
by applying Burkholder's inequality to the martingale difference sequence $(D_{n,k}^i-D_{n,k-1}^i)$, we deduce that for all $1\leqslant i \leqslant d$,
\begin{equation}\label{eq5-liang}
\mathbb{E} \left(\left|W_{n+1}^i-W_n^i\right|^{\alpha}\right)
\leqslant\mathbb{E}\left[\sup_{k}\left|D_{n,k}^i\right|^{\alpha}\right]
\leqslant B\mathbb{E}\left(\sum_{k=1}^{\infty}\left|D_{n,k}^i-D_{n,k-1}^i\right|^2\right)^{\alpha/2}.
\end{equation}
It is easy to see, from the definition of $D_{n,k}^i$ and the inequality $(x_1+\cdots+x_k)^p\leqslant k^{p-1}(x_1^p+\cdots +x_k^p)$ for $x_i\geqslant 0$ and $p>1$, that for all $1\leqslant i\leqslant d$,
\begin{eqnarray}\label{eq7-liang}
\left(\sum\limits_{k=1}^{\infty}\left|D_{n,k}^i-D_{n,k-1}^i\right|^2\right)^{\alpha/2}
&\leqslant & \left(\sum_{k=1}^{\infty}d \sum_{r=1}^{d} \frac{(u(r))^2}{(\rho^n u(i))^2}|W_{k,n}^r-1|^2\mathbf{1}_{\{Z_n^i(r)\geqslant k\}}\right)^{\alpha/2}\notag\\
&\leqslant& d^{\alpha-1}\sum_{r=1}^d  \left(\frac{u(r)}{\rho^n u(i)}\right)^{\alpha}\left(\sum_{k=1}^{Z_n^i(r)}|W_{k,n}^r-1|^2\right)^{\alpha/2};
\end{eqnarray}
and moreover, noting that $Z_n^i(r)$ are adapted to $\mathcal{F}_n$ and that $W_{k,n}^r$  are independent of $\mathcal{F}_n$ with the same distribution as $W_{1,0}^r$ for all $1\leqslant r\leqslant d$, it follows from Jensen's inequality that 
\begin{eqnarray}\label{eq4-liang}
&&\mathbb{E}\left[\left(\sum_{k=1}^{\infty}|D_{n,k}^i-D_{n,k-1}^i|^2\right)^{\alpha/2}|\mathcal{F}_n\right]\notag\\
&\leqslant& d^{\alpha-1}\sum_{r=1}^d  \left(\frac{u(r)}{\rho^n u(i)}\right)^{\alpha/2}\left(\frac{Z_n^i(r)u(r)}{\rho^n u(i)}\right)^{\alpha/2}\mathbb{E}\left[\left(\sum_{k=1}^{Z_n^i(r)}\frac{1}{Z_n^i(r)}|W_{k,n}^r-1|^2\right)^{\alpha/2}|\mathcal{F}_n\right]\notag\\
&\leqslant& d^{\alpha-1}\sum_{r=1}^d \left(\frac{u(r)}{\rho^n u(i)}\right)^{\alpha/2}(W_n^i)^{\alpha/2}\sum_{k=1}^{Z_n^i(r)}\frac{1}{Z_n^i(r)}\mathbb{E}|W_{k,n}^r-1|^{\alpha}\notag\\
&=& d^{\alpha-1}\sum_{r=1}^d \left(\frac{u(r)}{\rho^n u(i)}\right)^{\alpha/2}(W_n^i)^{\alpha/2}\sum_{k=1}^{Z_n^i(r)}\frac{1}{Z_n^i(r)}\mathbb{E}|W_{1,0}^r-1|^{\alpha}\notag\\
&\leqslant& d^{\alpha-1} \left(\max_{1\leqslant j\leqslant d} \mathbb{E}|W_{1,0}^j-1|^{\alpha}\right)(W_n^i)^{\alpha/2} \sum_{r=1}^d \left(\frac{u(r)}{\rho^n u(i)}\right)^{\alpha/2}\notag\\
&=& C (W_n^i)^{\alpha/2}\rho^{-\alpha n/2},
\end{eqnarray}
where $C=d^{\alpha-1}\left(\max_{1\leqslant j\leqslant d} \mathbb{E}|W_{1,0}^j-1|^{\alpha}\right) \sum_{r=1}^d \left(\frac{u(r)}{u(i)}\right)^{\alpha/2}$ is finite in view of   \eqref{eq3-liang} and  $u(r)\in (0,1)$ for all $1\leqslant r\leqslant d$.

It follows from \eqref{eq10}, \eqref{eq5-liang} and \eqref{eq4-liang} that
\begin{align}\label{eq11}
\sup_{n\geqslant0}\mathbb{E}|W^i_n-1|^\alpha&\leqslant C\left[\sum_{n=0}^{\infty}\left(\mathbb{E}(W_n^i)^{\alpha/2}\rho^{-\alpha n/2}\right)^{2/\alpha}\right]^{\alpha/2}\nonumber\\
&\leqslant C\sup_{n\geqslant0}\mathbb{E}(W_n^i)^{\alpha/2}\left[\sum_{n=0}^{\infty}\rho^{-n}\right]^{\alpha/2}.
\end{align}
As the assumption $\rho>1$ implies that $\sum_{n=0}^{\infty}\rho^{-n}=\frac{\rho}{\rho-1}<\infty$, the finiteness of $\sup_{n\geqslant 0}\mathbb{E}|W_n^i-1|^{\alpha}$ can be deduced by that of  $\sup_{n\geqslant 0}\mathbb{E}|W_n^i-1|^{\alpha/2}$.
In fact, for the positive integer $b$ such that $2^b\leqslant \alpha$, we can prove the following fact by induction that for all $1\leqslant i\leqslant d$,
\begin{align}\label{eq12}
\sup_{n\geqslant0}\mathbb{E}(W^i_n)^{\gamma}<\infty  \ \  \text{if} \ \ \mathbb{E}(W_{1}^i)^{\gamma}<\infty, \ \  \forall\gamma\in(2^{b-1},2^{b}].
\end{align}
Notice that $\mathbb{E}(W_1^i)^{\alpha}<\infty$ implies that $\mathbb{E}(W_1^i)^{\gamma}<\infty$ for all $1\leqslant \gamma\leqslant\alpha$.
If $b=1$, then $\gamma\in (1,2]$, using the result obtained for the case $\alpha\in(1,2]$ (cf. \eqref{eq6-liang}),  we have $\sup_{n\geqslant 0}\mathbb{E}|W_n^i-1|^{\gamma}<\infty$.
Assume that \eqref{eq12} holds true for some integer $b\geqslant1$.
Consider $\gamma\in(2^b,2^{b+1}]$. Since $\gamma/2\in (2^{b-1},2^b]$, \eqref{eq12} holds true for $\gamma/2$, that is, $\sup_{n\geqslant 0}\mathbb{E}(W_n^i)^{\gamma/2}<\infty$. Moreover, from \eqref{eq11}, we have $\sup_{n\geqslant 0}\mathbb{E}|W_n^i-1|^{\gamma}<\infty$.
Hence we get \eqref{eq12} holds true for all positive integers $b$ satisfying $2^b\leqslant \alpha$.
Furthermore, taking $\gamma=\alpha/2$, \eqref{eq12} implies that $\sup_{n\geqslant 0}\mathbb{E}(W_n^i)^{\alpha/2}<\infty$, and hence  $\sup_{n\geqslant 0}\mathbb{E}(W_n^i)^{\alpha}<\infty$  by \eqref{eq11},  which
completes the proof of the implication $(a)\Rightarrow(b)$.

The implications $(b)\Rightarrow(c)\Rightarrow(d)$ are evident. We now prove that $(d)\Rightarrow(a)$. Let $Z_{l, n, k}^r=(Z_{l, n, k}^r(1), \cdots, Z_{l, n, k}^r(d))$, where $Z_{l, n, k}^r(j)$ denotes the offspring of type $j$ at time $n+k$ of the $l$-th particle of type $r$ in the generation $n$. By iteration of (\ref{eq1}), it is easy to see that the MBP $(Z_n^i)$ satisfies the relation
$$
Z_{n+k}^i=\sum_{j=1}^d \sum_{l=1}^{Z_n^i(j)} Z_{l, n, k}^j, \quad n, k \geqslant 0 .
$$
We have, for all $n,k\geqslant0$ and $1\leqslant i\leqslant d,$
\begin{align}\label{eq14}
W_{n+k}^i&=\sum_{r=1}^d \frac{u(r)}{\rho^{n+k}u(i)} \sum_{j=1}^d \sum_{l=1}^{Z_n^i(j)} Z_{l, n, k}^j(r)\nonumber \\
& =\sum_{j=1}^d \frac{u(j)}{\rho^{n}u(i)} \sum_{l=1}^{Z_n^i(j)} \sum_{r=1}^d \frac{u(r) Z_{l, n, k}^j(r)}{\rho^{k}u(j)}\nonumber\\
& =\sum_{j=1}^d \frac{u(j)}{\rho^{n}u(i)} \sum_{l=1}^{Z_n^i(j)}W_{l, n, k}^j,
\end{align}
where
$$
W_{l, n, k}^j:=\frac{\langle Z_{l, n, k}^j, u\rangle}{\rho^k u(j)}.
$$
Since $$
Z_{l, n, k+1}^j=\sum_{r=1}^d \sum_{m=1}^{Z_{l, n, k}^j(r)} N_{m, n+k}^r,
$$
we get $
\mathbb{E}[Z_{l, n, k+1}^j \mid \mathcal{F}_{n+k}]=M^TZ_{l, n, k}^j
$, and then
\begin{align*}
\mathbb{E}[W_{l, n, k+1}^j\mid\mathcal{F}_{n+k}]=
\frac{\langle M^TZ_{l, n, k}^j, u\rangle}{\rho^{k+1} u(j)}=\frac{\langle Z_{l, n, k}^j, Mu\rangle}{\rho^{k+1} u(j)}
=\frac{\langle Z_{l, n, k}^j, \rho u\rangle}{\rho^{k+1} u(j)}=W_{l, n, k}^j.
\end{align*}
Hence $(W_{l, n, k}^j)_{k\geqslant0}$ is the fundamental martingale associated to the MBP $(Z_{l, n, k}^j)_{k\geqslant0}$, so that it converges a.s. to a random variable $W_{l, n, \infty}^j$ as $k\rightarrow\infty.$ Letting $k\rightarrow\infty$ in (\ref{eq14}), we get the following distributional equation on the limit variables $W^i$ for all $1\leqslant i\leqslant d$ and $n\geqslant0,$
\begin{align}\label{eq15}
W^i&=\sum_{j=1}^d \frac{u(j)}{\rho^{n}u(i)} \sum_{l=1}^{Z_n^i(j)}W_{l, n, \infty}^j
\end{align}
where the random variables $W_{l, n, \infty}^j$, indexed by $l\geqslant1,n \geqslant 0$ and $1 \leqslant j \leqslant d$, are independent of each other and independent of $Z_n^i$, each with the distribution
$$
\mathbb{P}\left(W_{l, n, \infty}^j \in \cdot\right)=\mathbb{P}\left(W^j \in \cdot\right) .
$$
Notice that  $\mathbb{E}W^i=c\in(0,\infty)$. Let $n=1$ in (\ref{eq15}), by Jensen's inequality,
\begin{align}
\mathbb{E}\left((W^i)^\alpha|\mathcal{F}_1\right)&=\mathbb{E}\left(\left(\sum_{r=1}^d \frac{u(r)}{\rho u(i)} \sum_{l=1}^{Z_1^i(r)} W_{l, 1,\infty}^r\right)^\alpha|\mathcal{F}_1\right)\nonumber\\
&\geqslant\left(\mathbb{E}\left(\sum_{r=1}^d \frac{u(r)}{\rho u(i)} \sum_{l=1}^{Z_1^i(r)} W_{l, 1,\infty}^r|\mathcal{F}_1\right)\right)^\alpha\nonumber\\
&=\left(\sum_{r=1}^d \frac{Z_1^i(r)u(r)}{\rho u(i)} \mathbb{E}W^r\right)^\alpha\nonumber\\
&=c^\alpha(W_1^i)^\alpha,
\end{align}
so that $$\mathbb{E}(W_1^i)^\alpha\leqslant c^{-\alpha}\mathbb{E}(W^i)^\alpha, \ \ \forall 1\leqslant i\leqslant d.$$
Since $\langle Z_1^i,u\rangle=\rho u(i) W_1^i$, and
\begin{eqnarray*}
    \langle Z_1^i,u\rangle^{\alpha}
    = \left(\sum_{r=1}^d Z_1^i(r) u(r)\right)^{\alpha}
    \geqslant  \sum_{r=1}^d (Z_1^i(r))^{\alpha} (u(r))^{\alpha},
\end{eqnarray*}
we have for all $1\leqslant i,j\leqslant d,$
\begin{equation*}
    \mathbb{E} (Z_1^i(j))^{\alpha}
    \leqslant \frac{1}{(u(j))^{\alpha}}\sum_{r=1}^d\mathbb{E} (Z_1^i(r))^{\alpha} (u(r))^{\alpha}\leqslant \frac{\mathbb{E} \langle Z_1^i, u\rangle^{\alpha}}{(u(j))^{\alpha}}=\frac{\mathbb{E} (W_1^i)^{\alpha}\rho^{\alpha} (u(i))^{\alpha}}{(u(j))^{\alpha}}<\infty.
\end{equation*}
\end{proof}
\section{Proof of Theorem \ref{Thm2}}
Owing to the double martingale structure described in Section 2, we can apply Burkholder-Davis-Gundy (BDG) inequalities (see e.g. \cite{CT2012}) to the martingale difference sequences $(W_{n+1}^i-W_n^i)$ and $(D_{n,k}^i-D_{n,k-1}^i)$, thereby studying the weighted moments of $W^i,1\leqslant i\leqslant d.$ For the reader's convenience, let us first recall these inequalities and state several consequences that will be used in our arguments. For a martingale sequence $\left\{\left(f_n, \mathcal{G}_n\right): n \geqslant 1\right\}$ defined on some probability space $(\Omega, \mathcal{G}, \mathbb{P})$, set $f_0=0, \mathcal{G}_0=\{\varnothing, \Omega\}, d_n=f_n-f_{n-1}$ for $n \geqslant 1$,
$$
\widehat{f}=\sup _{n \geqslant 1}\left|f_n\right| \quad \text { and } \quad \widehat{d}=\sup _{n \geqslant 1}\left|d_n\right|.
$$
\begin{lem}\label{lem3}
Let $\Phi$ be an increasing and continuous real function on $[0, \infty)$, with $\Phi(0)=0$ and $\Phi(2 \lambda) \leqslant c \Phi(\lambda)$ for some $c\in(0, \infty)$ and all $\lambda>0$.

(i) For every $\beta \in[1,2]$, there exists a constant $B=B_{c, \beta} \in(0, \infty)$ depending only on $c$ and $\beta$ such that for any martingale $\left\{\left(f_n, \mathcal{G}_n\right): n \geqslant 1\right\}$,
$$
\mathbb{E} \Phi\left(\widehat{f}\right) \leqslant B \mathbb{E} \Phi(s(\beta))+B \mathbb{E} \Phi\left(\widehat{d}\right), \quad \text { where } \quad s(\beta)=\left(\sum_{n=1}^{\infty} \mathbb{E}\left(\left|d_n\right|^\beta \mid \mathcal{G}_{n-1}\right)\right)^{1 / \beta},
$$
and
\begin{align}\label{eq16}
\mathbb{E} \Phi\left(\widehat{f}\right) \leqslant B \mathbb{E} \Phi(s(\beta))+B \sum_{n=1}^{\infty} \mathbb{E} \Phi\left(\left|d_n\right|\right).
\end{align}

(ii) If $\Phi$ is convex on $[0, \infty)$, then there exist constants $A=A_c \in(0, \infty)$ and $B=B_c \in(0, \infty)$ depending only on $c$ such that for any martingale $\left\{\left(f_n, \mathcal{G}_n\right): n \geqslant 1\right\}$,
$$
A \mathbb{E} \Phi(S) \leqslant \mathbb{E} \Phi\left(\widehat{f}\right) \leqslant B \mathbb{E} \Phi(S), \quad \text { where } \quad S=\left(\sum_{n=1}^{\infty} d_n^2\right)^{1 / 2}.
$$
Moreover, for any $\beta \in(0,2]$,
$$
\mathbb{E} \Phi\left(\widehat{f}\right) \leqslant B \mathbb{E} \Phi(S(\beta)), \quad \text { where } \quad S(\beta)=\left(\sum_{n=1}^{\infty}\left|d_n\right|^\beta\right)^{1 / \beta} .
$$

If, additionally, for some $\beta \in(0,2]$ the function $\Phi_{1 / \beta}(x):=\Phi\left(x^{1 / \beta}\right)$ is subadditive on $[0, \infty)$, then
$$
\mathbb{E} \Phi\left(\widehat{f}\right) \leqslant B \mathbb{E} \sum_{n=1}^{\infty} \Phi\left(\left|d_n\right|\right).
$$
\end{lem}
The following lemma concerns a smooth version of a regularly varying function, which allows us to apply convexity inequalities for martingales without altering the finiteness of the $\phi$-moment.
\begin{lem}\label{lem7}
(\cite{LL2013a}, Lemma 1.3) Let $\phi(x)=x^\alpha \ell(x)$, with $\alpha>1$, and $\ell(x)=\exp \left(\int_{a_0}^x \epsilon(u) d u / u\right)\left(x \geqslant a_0 \geqslant 0\right)$ with $\epsilon(x) \rightarrow 0(x \rightarrow \infty)$. Then, for each $\beta \in(1,2]$ with $\beta<\alpha$, there is a function $\phi_1(x) \geqslant 0$ such that
\begin{itemize}
\item[(i)] $\phi_1(x) \sim \phi(x)$ as $x \rightarrow \infty$;
\item[(ii)] $\phi_1(x)$ and $\phi_1\left(x^{1 / \beta}\right)$ are convex on $[0, \infty)$;
\item[(iii)] $\phi_1(x)=x^\alpha \ell_1(x)$, where $\ell_1(x)$ is slowly varying at $\infty$ and $\ell_1(x)>0$ for all $x \geqslant 0$.
\end{itemize}
\end{lem}
Applying BDG inequalities to the sequences $(W^i_{n+1}-W^i_n)$ and $(D_{n,k}^i-D_{n,k-1}^i)$, we get the following key lemma which will be used
in the proof of Theorem \ref{Thm2}.
\begin{lem}\label{lem4-liang}
Let $\phi$ be a convex and increasing function on $[0,\infty)$
with $\phi(0)=0$ and $\phi(2 x) \leqslant c \phi(x)$ for some constant $c \in(0, \infty)$ and all $x>0$. Let $\beta \in(1,2]$. If the function $x\mapsto\phi(x^{1/\beta})$ is convex, and $\max_{1\leqslant j\leqslant d}\mathbb{E}|W_{1, 0}^j-1|^{\beta}<\infty$, then
we have that for all $1\leqslant i\leqslant d,$
\begin{align}\label{eq17}
\mathbb{E} \phi\left(|\widehat{W}^i-1|\right)\leqslant C\sum_{n=0}^{\infty}
\frac{\mathbb{E}\phi\left((W_{n}^i)^{1/\beta}\right)}{A^i\left(\rho^n u(i)\right)^{\beta-1}}+C\sum_{n=0}^{\infty}\mathbb{E}\phi\left(\left(W_n^i\right)^{1/\beta}\frac{\max_{1\leqslant r\leqslant d}\left|W_{1,n}^r-1\right|}{\left(\rho^n u(i)\right)^{(\beta-1)/\beta}}\right),
\end{align}
where $A^i=\sum_{n=0}^{\infty}(\rho^n u(i))^{1-\beta}.$
\end{lem}
\begin{proof}
Using inequality (\ref{eq16}) to the martingale difference
$(W^i_{n+1}-W^i_{n})$, we have that for all $1\leqslant i\leqslant d,$
\begin{align}\label{eq19}
\mathbb{E}\phi\left(|\widehat{W}^i-1|\right)\leqslant B\mathbb{E}\phi\left[\left(\sum_{n=0}^{\infty}\mathbb{E}\left(|W_{n+1}^i-W_{n}^i|^\beta|\mathcal{F}_n\right)\right)^{1/\beta}\right]
+B\sum_{n=0}^{\infty}\mathbb{E}\phi\left(|W_{n+1}^i-W_{n}^i|\right).
\end{align}

We now estimate the first term on the right-hand side of the above inequality. Write $\sigma(\beta)=\max_{1\leqslant j\leqslant d}\mathbb{E}|W_{1, 0}^j-1|^{\beta}$ and $\phi_{1/\beta}(x)=\phi(x^{1/\beta})$.
Since $\phi(x)$ is increasing with $\phi(2x)\leqslant c\phi(x)$ for all $x>0$, for any constant $c_1>0$, we have $\phi\left(c_1 x\right) \leqslant c_2 \phi(x)$ for $c_2=c^{n_1}$ and any $x>0$, where $n_1$ is an integer such that $c_1 \leqslant 2^{n_1}$; moreover,
\begin{eqnarray*}
    \phi_{1/\beta}(c_1x)
    =\phi(c_1^{1/\beta}x^{1/\beta})
    \leqslant\phi(2^{n_1/\beta}x^{1/\beta})
    \leqslant c^{n_1}\phi_{1/\beta}(x).
\end{eqnarray*}
 Noting that $\sum_{n=0}^{\infty}\frac{1}{A^i(\rho^n u(i))^{\beta-1}}=1,$ it follows from  \eqref{eq7} and the convexity of the function $\phi_{1/\beta}(x)$
that for all $1\leqslant i\leqslant d,$
\begin{align}\label{eq20}
\mathbb{E}\phi\left[\left(\sum_{n=0}^{\infty}\mathbb{E}\left(|W_{n+1}^i-W_n^i|^{\beta}|\mathcal{F}_n\right)\right)^{1/\beta}\right]
&=\mathbb{E}\phi_{1/\beta}\left[\sum_{n=0}^{\infty}\mathbb{E}\left(|W_{n+1}^i-W_n^i|^{\beta}|\mathcal{F}_n\right)\right]\nonumber\\
&\leqslant\mathbb{E}\phi_{1/\beta}\left[\sum_{n=0}^{\infty}\frac{1}{\left(\rho^n u(i)\right)^{\beta-1}}B_{\beta}^{\beta}W_{n}^i\sigma(\beta)\right]\nonumber\\
&\leqslant\mathbb{E}\sum_{n=0}^{\infty}\frac{1}{A^i\left(\rho^n u(i)\right)^{\beta-1}}\phi_{1/\beta}\left(A^iB_{\beta}^{\beta}W_{n}^i\sigma(\beta)\right)\nonumber\\
&\leqslant C\sum_{n=0}^{\infty}\frac{\mathbb{E}\phi\left((W_{n}^i)^{1/\beta}\right)}{A^i\left(\rho^n u(i)\right)^{\beta-1}}.
\end{align}

 For the second term on the right-hand side of (\ref{eq19}), applying Lemma \ref{lem3} (ii) to the martingale difference
$(D_{n,k}^i-D_{n,k-1}^i)$ and using the monotonicity of $\phi$ , we obtain that for all $1\leqslant i\leqslant d,$
\begin{eqnarray}\label{eq8-liang}
\mathbb{E}\left[\phi(|W_{n+1}^i-W_{n}^i|)|\mathcal{G}_{n,k-1}\right]
&\leqslant&\mathbb{E}\left[\phi(\sup_{k}|D_{n,k}^i|)|\mathcal{G}_{n,k-1}\right]\notag\\
&\leqslant &B\mathbb{E}\left[\phi\left(\left(\sum_{k=1}^{\infty}
\left|D_{n,k}^i-D_{n,k-1}^i\right|^\beta\right)^{1/\beta}\right)|\mathcal{G}_{n,k-1}\right].
\end{eqnarray}
Similarly to the argument of \eqref{eq7-liang}, as the function $x\mapsto x^{\beta}$ ($\beta\in(1,2]$) is convex, we have that for all $1\leqslant i\leqslant d$,
\begin{eqnarray}\label{XEQ4}
\sum\limits_{k=1}^{\infty}\left|D_{n,k}^i-D_{n,k-1}^i\right|^{\beta}
&\leqslant & \sum_{k=1}^{\infty}  \sum_{r=1}^{d} \frac{u(r)}{(\rho^n u(i))^{\beta}}|W_{k,n}^r-1|^{\beta}\mathbf{1}_{\{Z_n^i(r)\geqslant k\}}\nonumber\\
&\leqslant& W_n^i \max_{1\leqslant j\leqslant d}\sum_{k=1}^{Z_n^i(j)}\frac{1}{Z_n^i(j)}\frac{|W_{k,n}^j-1|^{\beta}}{(\rho^n u(i))^{\beta-1}}.
\end{eqnarray}
Moreover, by the monotonicity and convexity of $\phi_{1/\beta}$,
\begin{eqnarray}\label{eq9-liang}
\mathbb{E}\left[\phi_{1/\beta}\left(\sum_{k=1}^{\infty}|D_{n,k}^i-D_{n,k-1}^i|^{\beta}\right)|\mathcal{F}_{n}\right]
&\leqslant& \mathbb{E}\left[\phi_{1/\beta}\left(W_n^i \max_{1\leqslant j\leqslant d}\sum_{k=1}^{Z_n^i(j)}\frac{1}{Z_n^i(j)}\frac{|W_{k,n}^j-1|^{\beta}}{(\rho^n u(i))^{\beta-1}}\right)|\mathcal{F}_{n}\right]\notag\\
&\leqslant& \sum_{j=1}^d \mathbb{E}\left[\phi_{1/\beta}\left( \sum_{k=1}^{Z_n^i(j)}\frac{1}{Z_n^i(j)}W_n^i\frac{|W_{k,n}^j-1|^{\beta}}{(\rho^n u(i))^{\beta-1}}\right)|\mathcal{F}_{n}\right]\notag\\
&\leqslant& \sum_{j=1}^d \sum_{k=1}^{Z_n^i(j)}\frac{1}{Z_n^i(j)} \mathbb{E}\left[\phi_{1/\beta}\left(W_n^i \frac{|W_{k,n}^j-1|^{\beta}}{(\rho^n u(i))^{\beta-1}}\right)|\mathcal{F}_{n}\right]\notag\\
&=&\sum_{j=1}^d \mathbb{E}\left[\phi_{1/\beta}\left(W_n^i \frac{|W_{1,n}^j-1|^{\beta}}{(\rho^n u(i))^{\beta-1}}\right)|\mathcal{F}_{n}\right].
\end{eqnarray}
Taking expectations on both sides of \eqref{eq8-liang}, it follows from \eqref{eq9-liang}
that for all $1\leqslant i\leqslant d,$
\begin{align}\label{eq18}
\mathbb{E}\phi\left(|W_{n+1}^i-W_{n}^i|\right)&\leqslant B\sum_{j=1}^d\mathbb{E}\left[\phi
\left((W_{n}^i)^{1/\beta}\frac{|W_{1,n}^j-1|}{(\rho^n u(i))^{(\beta-1)/\beta}}\right)\right]\nonumber\\
&\leqslant Bd\mathbb{E}\phi\left(\left(W_n^i\right)^{1/\beta}\frac{\max_{1\leqslant l\leqslant d}\left|W_{1,n}^l-1\right|}{\left(\rho^n u(i)\right)^{(\beta-1)/\beta}}\right).
\end{align}
Combining (\ref{eq19}), (\ref{eq20}) and \eqref{eq18}, we complete the proof of (\ref{eq17}).
\end{proof}

\begin{proof}[Proof of Theorem \ref{Thm2}]  Let $\beta \in(1,2]$ with $\beta<\alpha$. Write $\phi(x)=x^\alpha \ell(x)$. By Lemma \ref{lem7}, we can assume that both $\phi(x)$ and $\phi(x^{1 / \beta})$ are increasing and convex on $[0, \infty)$ with $\phi(2x)\leqslant c\phi(x)$ and
$\ell(x)>0$ for some $c>0$ and all $x \geqslant 0$.

By the convexity of $\phi(x)$ ($\sum_{j=1}^d u(j)=1$), we have that for all $1\leqslant i\leqslant d,$
\begin{align*}
\mathbb{E}\phi(W^i_1)&=\mathbb{E}\phi\left(\sum_{j=1}^d\frac{Z^i_1(j)u(j)}{\rho u(i)}\right)
\leqslant\sum_{j=1}^du(j)\mathbb{E}\phi\left(\frac{Z^i_1(j)}{\rho u(i)}\right)
\leqslant C\sum_{j=1}^du(j)\mathbb{E}\phi\left(Z^i_1(j)\right),
\end{align*}
where the last inequality is according to the condition $\phi(2x)\leqslant c\phi(x)$.
This completes the proof of (a)$\Rightarrow$(b).

We now prove that (b)$\Rightarrow$(c). By Lemma \ref{lem4-liang}, we have for all $1\leqslant i\leqslant d,$
\begin{align}\label{eq10-liang}
\mathbb{E}\phi\left(|\widehat{W}^i-1|\right)\leqslant C\sum_{n=0}^{\infty}I^i_{n}+C\sum_{n=0}^{\infty}J^i_n,
\end{align}
where $$I^i_n=\frac{\mathbb{E}\phi\left((W_{n}^i)^{1/\beta}\right)}{A^i\left(\rho^n u(i)\right)^{\beta-1}},$$
$$J^i_n=\mathbb{E}\phi\left(\left(W_n^i\right)^{1/\beta}\frac{\max_{1\leqslant r\leqslant d}\left|W_{1,n}^r-1\right|}{\left(\rho^n u(i)\right)^{(\beta-1)/\beta}}\right).$$

We first show the summability of the terms $(I_n^i)_{n\geqslant 0}$ for all $1\leqslant i\leqslant d$, that is,
\begin{equation}\label{eq11-liang}
\sum_{n=0}^{\infty}I^i_n<\infty\,.
\end{equation}
Since $\ell$ is bounded away from 0 and $\infty$ on any compact subset of $[0, \infty)$, by Potter's theorem (see \cite{BGT1987}), for any
$\delta>0$, there exists $C=C(\ell, \delta)>1$ such that $\ell(x) \leqslant C \max (x^\delta, x^{-\delta}) \leqslant C(x^\delta+x^{-\delta})$ for all $x>0$. Hence
\begin{align}\label{eq21}
I^i_n&=\frac{\mathbb{E}(W_{n}^i)^{\alpha/\beta}\ell\left((W_{n}^i\right)^{1/\beta})}{A^i\left(\rho^n u(i)\right)^{\beta-1}}\leqslant C\frac{\mathbb{E}(W_{n}^i)^{\frac{\alpha+\delta}{\beta}}}{A^i\left(\rho^n u(i)\right)^{\beta-1}}+C\frac{\mathbb{E}(W_{n}^i)^{\frac{\alpha-\delta}{\beta}}}{A^i\left(\rho^n u(i)\right)^{\beta-1}}.
\end{align}
By Potter's theorem again, for the same $\delta>0$, there exists $C=C(\ell, \delta)>0$ such that $\ell(x) \geqslant C x^{-\delta}$ for all $x \geqslant 1$, so that $x^{\alpha-\delta} \leqslant 1+\frac{1}{C} x^\alpha \ell(x)=1+\frac{1}{C} \phi(x)$ for all $x \geqslant 0$ when $\delta<\alpha$.
Therefore
$$
\mathbb{E}(W^i_1)^{\alpha-\delta} \leqslant 1+\frac{1}{C} \mathbb{E} \phi\left(W^i_1\right)<\infty.
$$
Choose $0<\delta\leqslant\alpha\frac{\beta-1}{\beta+1}$, by H\"{o}lder's inequality, we have that for all $1\leqslant i\leqslant d,$
$$\left[\mathbb{E}(W_1^i)^{\frac{\alpha+\delta}{\beta}}\right]^{\frac{\beta}{\alpha+\delta}}
\leqslant\left[\mathbb{E}(W_1^i)^{\alpha-\delta}\right]^{\frac{1}{\alpha-\delta}}<\infty.$$
Moreover, it follows from Theorem \ref{Thm1} that
\begin{align}\label{XEQ5}
\sup_{n\geqslant0}\mathbb{E}(W_{n}^i)^{\frac{\alpha+\delta}{\beta}}<\infty\,, \ \ \forall1\leqslant i\leqslant d.
\end{align}
Similarly, we can also obtain
$$
\sup_{n\geqslant0}\mathbb{E}(W_{n}^i)^{\frac{\alpha-\delta}{\beta}}<\infty\,, \ \ \forall1\leqslant i\leqslant d.
$$
Since $A^i\in(0,\infty)$, $\rho>1$ and $u(i)\in(0,1)$, we obtain from \eqref{eq21}
that
\begin{align*}
\sum_{n=0}^\infty I^i_n\leqslant C\sum_{n=0}^\infty(\rho^nu(i))^{1-\beta}<\infty\,, \ \ \forall1\leqslant i\leqslant d.
\end{align*}

It remains to prove that
\begin{align}\label{eq22}
\sum_{n=0}^\infty J^i_n<\infty\,, \ \ \forall1\leqslant i\leqslant d.
\end{align}
Again by Potter's theorem, for $\delta>0$, there exists $C=C(\ell, \delta)>0$ such that $\ell(x y) \leqslant C \ell(x) \max \left\{y^\delta, y^{-\delta}\right\} \leqslant C \ell(x)\left(y^\delta+y^{-\delta}\right)$ for all $x>0$, $y>0$. Recall that the random variables $W_{1, n}^r$ are independent of $Z_n^i$, we have
\begin{eqnarray*}
J_n^i
&=&\mathbb{E}\left(\left(W_n^i\right)^{\alpha/\beta}\frac{\left(\max_{1\leqslant r\leqslant d}\left|W_{1,n}^r-1\right|\right)^\alpha}{\left(\rho^n u(i)\right)^{\alpha(\beta-1)/\beta}}
\ell\left(\left(W_n^i\right)^{1/\beta}\frac{\max_{1\leqslant r\leqslant d}\left|W_{1,n}^r-1\right|}{\left(\rho^n u(i)\right)^{(\beta-1)/\beta}}\right)\right)\\
&\leqslant& C\mathbb{E}\left(W_n^i\right)^{\alpha/\beta}\frac{\left(\max_{1\leqslant r\leqslant d}\left|W_{1,n}^r-1\right|\right)^\alpha}{\left(\rho^n u(i)\right)^{\alpha(\beta-1)/\beta}}
\ell\left(\max_{1\leqslant r\leqslant d}\left|W_{1,n}^r-1\right|\right)\\
& & \qquad \cdot\left[\left(\frac{(W_n^i)^{1/\beta}}{(\rho^n u(i))^{(\beta-1)/\beta}}\right)^{\delta}\right.
 +\left.\left(\frac{(W_n^i)^{1/\beta}}{(\rho^n u(i))^{(\beta-1)/\beta}}\right)^{-\delta}\right]\\
&=& C\left[\left(\rho^n u(i)\right)^{\frac{-(\alpha+\delta)(\beta-1)}{\beta}}\mathbb{E}\left(W_n^i\right)^{(\alpha+\delta)/\beta} + \left(\rho^n u(i)\right)^{\frac{-(\alpha-\delta)(\beta-1)}{\beta}}\mathbb{E}\left(W_n^i\right)^{(\alpha-\delta)/\beta} \right] \mathbb{E}\phi\left(\max_{1\leqslant r\leqslant d}\left|W_{1,0}^r-1\right|\right).
\end{eqnarray*}
Since $\mathbb{E} \phi(W^i_1)<\infty$ ($\forall1\leqslant i\leqslant d$)  implies that $\mathbb{E}\phi(\max_{1\leqslant r\leqslant d}|W_{1,0}^r-1|)<\infty$, combining with the facts that $\sup_{n\geqslant0}\mathbb{E}(W_{n}^i)^{\frac{\alpha+\delta}{\beta}}<\infty$ and $\sup_{n\geqslant0}\mathbb{E}(W_{n}^i)^{\frac{\alpha-\delta}{\beta}}<\infty$, we see that
\begin{equation*}
    \sum_{n=0}^{\infty} J_n^i \leqslant C\left[\sum_{n=0}^{\infty}(\rho^nu(i))^{\frac{-(\alpha+\delta)(\beta-1)}{\beta}}+\sum_{n=0}^{\infty}(\rho^n u(i))^{\frac{-(\alpha-\delta)(\beta-1)}{\beta}}\right]<\infty\,,
\end{equation*}
where the last finiteness follows from $\rho>1$, $u(i)>0$ and ${\frac{(\alpha+\delta)(\beta-1)}{\beta}}>{\frac{(\alpha-\delta)(\beta-1)}{\beta}}>0$.

It follows from \eqref{eq10-liang}, \eqref{eq11-liang} and \eqref{eq22} that
$\mathbb{E}\phi(|\widehat{W}^i-1|)<\infty$, which is equivalent to $\mathbb{E}\phi(\widehat{W}^i)<\infty$ for all $1 \leqslant i\leqslant d$.
This completes the proof of (b)$\Rightarrow$(c).

We then show that $(c) \Rightarrow (d)$. Obviously, for all $1\leqslant i\leqslant d,$
$$
\mathbb{E} \phi(W^i) \leqslant \mathbb{E} \phi(\widehat{W}^i)<\infty ;
$$
by the dominated convergence theorem, we see that for all $1\leqslant i\leqslant d,$
$$
\mathbb{E}W^i=\lim _{n \rightarrow \infty} \mathbb{E} W_n^i=1.
$$

We finally show that $(d)\Rightarrow(a)$. The arguments are similar to those
in Theorem \ref{Thm1}.
Since the function $x\mapsto \phi(x)$ is convex and takes the same role as $x\mapsto x^{\alpha}$,
\begin{equation*}
    \mathbb{E}\phi(W_1^i)\leqslant \mathbb{E}\phi(W^i)<\infty\,.
\end{equation*}
Using Potter's theorem, together with
the monotonicity of $\phi$, we obtain that
\begin{align*}
\mathbb{E} \phi(Z_1^i(j))&=\mathbb{E} (Z_1^i(j))^{\alpha}\ell(Z_1^i(j))\\
&= \left(\frac{\rho u(i)}{u(j)}\right)^{\alpha}\mathbb{E} \left(\frac{Z_1^i(j) u(j)}{\rho u(i)}\right)^{\alpha}\ell\left(\frac{Z_1^i(j) u(j)}{\rho u(i)} \cdot \frac{\rho u(i)}{u(j)}\right) \\
&\leqslant C \left(\frac{\rho u(i)}{u(j)}\right)^{\alpha} \left[\left(\frac{\rho u(i)}{u(j)}\right)^{\delta}+\left(\frac{\rho u(i)}{u(j)}\right)^{-\delta}\right]\mathbb{E} \left(\frac{Z_1^i(j) u(j)}{\rho u(i)}\right)^{\alpha}\ell\left(\frac{Z_1^i(j) u(j)}{\rho u(i)}\right) \\
&=C \left[\left(\frac{\rho u(i)}{u(j)}\right)^{\alpha+\delta}+\left(\frac{\rho u(i)}{u(j)}\right)^{\alpha-\delta}\right] \mathbb{E} \phi\left(\frac{Z_1^i(j) u(j)}{\rho u(i)}\right) \\
&\leqslant C \left[\left(\frac{\rho u(i)}{u(j)}\right)^{\alpha+\delta}+\left(\frac{\rho u(i)}{u(j)}\right)^{\alpha-\delta}\right] \mathbb{E}\phi(W_1^i)<\infty\,.
\end{align*}
\end{proof}

\section{Proof of Theorem \ref{Thm3}}

Before giving the proof of Theorem \ref{Thm3}, we first present a lemma which plays a similar role to Lemma \ref{lem4-liang} for the case $\alpha>1$:
\begin{lem}\label{Lm5}
    Let $\phi:[0,\infty) \mapsto [0,\infty)$ be convex and increasing with $\phi(0)=0$, $\phi(2x)\leqslant c\phi(x)$ for some constant $c\in(0,\infty)$ and all $x\geqslant 0$.
    If the function $x\mapsto \phi(x^{1/2})$
    is concave, then
    \begin{equation*}
        \mathbb{E} \phi(\widehat{W}^i-1)\leqslant C\sum_{n=0}^{\infty} \rho^n u(i) \mathbb{E}\phi\left(\max_{1\leqslant j\leqslant d}\frac{ |W_{1,0}^j-1|}{\rho^n u(i)}\right).
    \end{equation*}
\end{lem}

\begin{proof}
Note that the function $\phi_{1/2}(x):=\phi(x^{1/2})$ is concave, which implies its sub-additivity. Hence applying Lemma \ref{lem3} (ii) (with $\beta=2$)
 to the martingale difference $(W_{n+1}^i-W_n^i)$,
we see that
\begin{equation}\label{eq12-liang}
   \mathbb{E} \phi(\widehat{W}^i-1)\leqslant B \mathbb{E}\sum_{n=0}^{\infty} \phi(|W_{n+1}^i-W_n^i|).
\end{equation}

 Similarly to the argument in  \eqref{eq7-liang}, the sub-additivity of $\phi_{1/2}(x)$ and the convexity of $\phi(x)$ yield that
\begin{eqnarray*}
\phi_{1/2}\left(\sum_{k=1}^{\infty}|D_{n,k}^i-D_{n,k-1}^i|^2\right)
&\leqslant & \sum_{k=1}^{\infty} \phi_{1/2}(|D_{n,k}^i-D_{n,k-1}^i|^2)\\
&\leqslant& \sum_{k=1}^{\infty}\sum_{r=1}^d u(r)\phi\left(\frac{|W_{k,n}^r-1|}{\rho^n u(i)}\mathbf{1}_{\{Z_n^i(r)\geqslant k\}}\right)\\
&=& \sum_{r=1}^d\sum_{k=1}^{Z_n^i(r)} u(r)\phi\left(\frac{|W_{k,n}^r-1|}{\rho^n u(i)}\right).
\end{eqnarray*}
 Recall that
$W_{n+1}^i-W_n^i$ is the sum of martingale difference sequence $(D_{n,k}^i-D_{n,k-1}^i)$,
 using  Lemma \ref{lem3} (ii) again to the martingale difference $(D_{n,k}^i-D_{n,k-1}^i)$, together with the fact that $(W_{k,n}^r)$ are i.i.d. and are independent of $\mathcal{F}_n$,
we obtain that for all $1\leqslant i\leqslant d$,
\begin{align*}
    \mathbb{E}\left[\phi(|W_{n+1}^i-W_{n}^i|)|\mathcal{F}_n\right]
    &\leqslant B\mathbb{E}\left[\phi_{1/2}\left(\sum_{k=1}^{\infty}|D_{n,k}^i-D_{n,k-1}^i|^2 \right) |\mathcal{F}_n\right]\\
    & \leqslant B \mathbb{E} \left[\sum_{r=1}^d\sum_{k=1}^{Z_n^i(r)} u(r)\phi\left(\frac{|W_{1,n}^r-1|}{\rho^n u(i)}\right)|\mathcal{F}_n\right]\\
    &\leqslant B\sum_{r=1}^d Z_{n}^i(r) u(r)\mathbb{E}\phi\left(\max_{1\leqslant j\leqslant d}\frac{ |W_{1,n}^j-1|}{\rho^n u(i)}\right)\\
    &= B \rho^n u(i)W_n^i \mathbb{E} \phi\left(\max_{1\leqslant j\leqslant d}\frac{ |W_{1,n}^j-1|}{\rho^n u(i)}\right).
\end{align*}
Taking expectation on both sides,
\begin{align*}
    \mathbb{E} \phi\left(|W_{n+1}^i-W_n^i|\right)
    &\leqslant B \rho^n u(i) \mathbb{E} W_n^i \mathbb{E}\phi\left(\max_{1\leqslant j\leqslant d}\frac{ |W_{1,n}^j-1|}{\rho^n u(i)}\right)\\
    &= B\rho^n u(i) \mathbb{E} \phi\left(\max_{1\leqslant j\leqslant d}\frac{|W_{1,0}^j-1|}{\rho^n u(i)}\right),
\end{align*}
and moreover, combining with \eqref{eq12-liang},
\begin{equation*}
    \mathbb{E}\phi(\widehat{W}^i-1)
    \leqslant C\sum_{n=0}^{\infty} \rho^n u(i)\mathbb{E} \phi\left(\max_{1\leqslant j\leqslant d}\frac{ |W_{1,0}^j-1|}{\rho^n u(i)}\right).
\end{equation*}
\end{proof}

\begin{proof}[Proof of Theorem \ref{Thm3}]
Let $\phi(x)=x\ell(x)$ and $\tilde{\phi}(x)=x \tilde{\ell}(x)$.
Without loss of generality, we can assume that $\phi(x)=x\ell(x)$ is convex with $\phi(0)=0$, $\phi(2x)\leqslant c\phi(x)$ for some constant $c\in (0,\infty)$
 and all $x>0$,
$\phi_{1/2}(x)$ is concave on $(0,\infty)$, and moreover, $\ell(x)$ is linear in a neighborhood of $0$.

Notice that $\ell(x)$ is increasing on $[0,\infty)$, it follows from Lemma \ref{Lm5} and Potter's theorem that
\begin{align*}
    \mathbb{E}\phi(\widehat{W}^i-1)
    &\leqslant C\sum_{n=0}^{\infty} \rho^n u(i)\mathbb{E} \phi\left(\frac{\max |W_{1,0}^r-1|}{\rho^n u(i)}\right)\\
    &=C \sum_{n=0}^{\infty} \rho^n u(i) \mathbb{E} \frac{\max|W_{1,0}^r-1|}{\rho^n u(i)}\ell\left(\frac{\max|W_{1,0}^r-1|}{\rho^n u(i)}\right)\\
    &\leqslant C\sum_{n=0}^{\infty}\mathbb{E} \max|W_{1,0}^r-1|\ell\left(\rho^{-n}\max |W_{1,0}^r-1|\right)\\
    &\leqslant C \mathbb{E} \max|W_{1,0}^r-1|\ell(\max|W_{1,0}^r-1|)+  C \sum_{n=1}^{\infty}\mathbb{E} \max|W_{1,0}^r-1|\int_{\rho^{-n}\max|W_{1,0}^r-1|}^{\rho^{1-n}\max |W_{1,0}^r-1|}\frac{\ell(t)}{t}dt\\
    &= C\mathbb{E} \phi( \max|W_{1,0}^r-1|)+C \mathbb{E} \max|W_{1,0}^r-1|\left(\int_{0}^1\frac{\ell(t)}{t}dt+\int_1^{\max|W_{1,0}^r-1|} \frac{\ell(t)}{t}dt\right)\\
    &\leqslant C\mathbb{E} \phi(\max |W_{1,0}^r-1|)+ C\mathbb{E} \max|W_{1,0}^r-1|+C \mathbb{E} \tilde{\phi}(\max |W_{1,0}^r-1|).
\end{align*}
The finiteness of $\int_0^1 \frac{\ell(t)}{t}dt$ is due to the fact that $\ell(x)$ is linear in a neighborhood of $0$.

In the following, we will show the finiteness of the term $\mathbb{E}\tilde{\phi}(\max |W_{1,0}^r-1|)$.
Notice that
$\ell(x)$ is increasing with $\ell(2x)\leqslant c\ell(x)$ for all $x>0$, it follows from the definition of $\tilde{\ell}(x)$ that when $x\geqslant2$,
\begin{align*}
    \tilde{\ell}(x)
    &= \int_{1}^x\frac{\ell(t)}{t}dt
    \geqslant \int_{x/2}^x \frac{\ell(t)}{t}dt
    \geqslant \frac{\ell(x/2)}{x}\cdot \frac{x}{2}
    = \frac{1}{2}\ell(x/2)
    \geqslant \frac{1}{2c}\ell(x),
\end{align*}
and moreover, $\phi(x)\leqslant  2c \tilde{\phi}(x)$ for $x\geqslant 2$.
Again by the monotonicity of  $\ell(x)$, when $x\geqslant 1$,
\begin{equation*}
    \tilde{\ell}(2x)=\int_1^{2x}\frac{\ell(t)}{t}dt =\int_1^x\frac{\ell(t)}{t}dt +\int_x^{2x}\frac{\ell(t)}{t}dt \leqslant \tilde{\ell}(x)+\ln 2 \ell(2x),
\end{equation*}
and hence,
\begin{equation*}
    \tilde{\phi}(2x)\leqslant 2\tilde{\phi}(x)+\ln 2 \phi(2x).
\end{equation*}
Since $|W_{1,0}^r-1|\leqslant W_{1,0}^r+1 \leqslant 2\max\{W_{1,0}^r,1\}$,
 together with the facts that
$\tilde{\phi}(x)$ is increasing on $[0,\infty)$ and $\phi(2 x)\leqslant c\phi(x)$, by the assumption $\mathbb{E}\tilde{\phi}(W_{1,0}^r)=\mathbb{E}\tilde{\phi}(W_1^r)<\infty$, we have for all $1\leqslant r\leqslant d$,
\begin{align*}
    \mathbb{E}\tilde{\phi}(|W_{1,0}^r-1|)
    &\leqslant \mathbb{E} \tilde{\phi}(2\max\{W_{1,0}^r,1\})\\
    &\leqslant 2\mathbb{E}\tilde{\phi}(\max\{W_{1,0}^r,1\})+\ln 2 \mathbb{E} \phi(2\max\{W_{1,0}^r,1\})\\
    &\leqslant 2\mathbb{E}\tilde{\phi}(\max\{W_{1,0}^r,1\})+c\ln 2 \mathbb{E}\phi(\max\{W_{1,0}^r,1\})\\
    &\leqslant  2\mathbb{E}\tilde{\phi}(W_{1,0}^r)+c\ln 2 \mathbb{E}\phi(W_{1,0}^r)\mathbf{1}_{\{W_{1,0}^r>2\}}+2\tilde{\phi}(1)+c\ln 2 \phi(1)+c\ln 2\phi(2)\\
    &\leqslant  (2+2c^2\ln 2)\mathbb{E} \tilde{\phi}(W_{1,0}^r)+(2\tilde{\phi}(1)+c\ln 2\phi(1)+c\ln 2\phi(2))\\
    &<\infty\,;
\end{align*}
 moreover,
\begin{equation}\label{eq13-liang}
    \mathbb{E}\tilde{\phi}(\max|W_{1,0}^r-1|)
    \leqslant \sum_{r=1}^d\mathbb{E}\tilde{\phi}(|W_{1,0}^r-1|)<\infty \,.
\end{equation}

On the other hand, the terms $\mathbb{E} \phi(\max|W_{1,0}^r-1|)$ and $\mathbb{E} \max |W_{1,0}^r-1|$ are finite as
\begin{align*}
    \mathbb{E} \phi(\max |W_{1,0}^r-1|)
    &\leqslant \sum_{r=1}^d\mathbb{E} \phi(|W_{1,0}^r-1|)\\
    &= \sum_{r=1}^d \left(\mathbb{E}\phi(|W_{1,0}^r-1|)\mathbf{1}_{\{|W_{1,0}^r-1|>2\}} +\mathbb{E}\phi(|W_{1,0}^r-1|)\mathbf{1}_{\{|W_{1,0}^r-1|\leqslant2\}}\right)\\
    &\leqslant 2c \sum_{r=1}^d\mathbb{E} \tilde{\phi}(|W_{1,0}^r-1|)\mathbf{1}_{\{|W_{1,0}^r-1|>2\}} + d\phi(2)\\
    &\leqslant 2c \sum_{r=1}^d\mathbb{E} \tilde{\phi}(|W_{1,0}^r-1|)+d\phi(2)<\infty\,,
\end{align*}
and
\begin{align*}
    \mathbb{E} \max |W_{1,0}^r-1|
    &\leqslant \sum_{r=1}^d\mathbb{E} |W_{1,0}^r-1|\\
    &= \sum_{r=1}^d\left(\mathbb{E} |W_{1,0}^r-1|\mathbf{1}_{\{|W_{1,0}^r-1|\geqslant 1\}}+\mathbb{E}|W_{1,0}^r-1|\mathbf{1}_{\{|W_{1,0}^r-1|< 1\}}\right)\\
    &\leqslant   \sum_{r=1}^d\mathbb{E} |W_{1,0}^r-1|\mathbf{1}_{\{|W_{1,0}^r-1|\geqslant 1\}}+d\\
    &\leqslant   \sum_{r=1}^d\frac{1}{\ell(1)}\mathbb{E}\phi(|W_{1,0}^r-1|)+d<\infty\,.
\end{align*}

Therefore, we show that
\begin{equation*}
   \mathbb{E}\phi(\widehat{W}^i-1)<\infty\,,
\end{equation*}
which is equivalent to $\mathbb{E}\phi(\widehat{W}^i)<\infty$.
The finiteness of $\mathbb{E}\phi(W^i)$ is obvious.
\end{proof}

\section{Proof of Theorem \ref{Thm6}}
In this section, we consider the upper estimate of convergence rates in the $\phi$-moment sense, based on Theorem 1.1 and the BDG inequalities to the martingale difference sequence
$\{\delta^n(D_{n,k}^i-D_{n,k-1}^i)\}$.

\begin{proof}[Proof of Theorem \ref{Thm6}]
$(i)$ Let $\beta\in(1,2]$ and $\beta<\alpha.$ By a standard modification of regularly varying functions, we may replace $\phi$ by an equivalent function, still denoted by $\phi$, such that both $\phi(x)$ and $\phi(x^{1 / \beta})$ are increasing and convex on $[0, \infty)$, and $\phi(2x)\leqslant c\phi(x)$  for some $c>0$ and all $x\geqslant 0$. Applying Lemma \ref{lem3} (ii) to the martingale difference sequence
$\{\delta^n(D_{n,k}^i-D_{n,k-1}^i)\}$
with respect to the filtration $\{\mathcal{G}_{n,k-1}\}$,
and using the monotonicity of $\phi$ , we obtain that for all $1\leqslant i\leqslant d,$
\begin{eqnarray}\label{XEQ12}
\mathbb{E}\left[\phi(\delta^n|W_{n+1}^i-W_{n}^i|)|\mathcal{F}_{n}\right]
&\leqslant&\mathbb{E}\left[\phi(\delta^n\sup_{k}|D_{n,k}^i|)|\mathcal{F}_{n} \right]\notag\\
&\leqslant &B\mathbb{E}\left[\phi_{1/\beta}\left(\sum_{k=1}^{\infty}
\delta^{\beta n}\left|D_{n,k}^i-D_{n,k-1}^i\right|^\beta\right)|\mathcal{F}_{n}\right].
\end{eqnarray}
By an argument similar to that used in \eqref{eq9-liang}, it follows that
\begin{align}\label{10eq2}
\mathbb{E}\phi\left(\delta^n|W_{n+1}^i-W_{n}^i|\right)\leqslant Bd\mathbb{E}\phi\left(\delta^n\left(W_n^i\right)^{1/\beta}\frac{\max_{1\leqslant l\leqslant d}\left|W_{1,n}^l-1\right|}{\left(\rho^n u(i)\right)^{(\beta-1)/\beta}}\right).
\end{align}
By Potter's theorem, for $\varepsilon>0$, there exists $C=C(\ell, \varepsilon)>0$ such that $\ell(x y) \leqslant C \ell(x) \max \left\{y^\varepsilon, y^{-\varepsilon}\right\} \leqslant C \ell(x)\left(y^\varepsilon+y^{-\varepsilon}\right)$ for all $x>0$, $y>0$. Using the fact that $W_{1, n}^r$ is independent of {\color{blue} $W_n^i$}, we have
\begin{eqnarray}\label{10eq3}
&&\mathbb{E}\phi\left(\delta^n\left(W_n^i\right)^{1/\beta}\frac{\max_{1\leqslant l\leqslant d}\left|W_{1,n}^l-1\right|}{\left(\rho^n u(i)\right)^{(\beta-1)/\beta}}\right)\nonumber\\
&\leqslant& C\left[\left(\frac{\delta}{\rho^{1-1/\beta}}\right)^{(\alpha+\varepsilon)n}\mathbb{E}\left(W_n^i\right)^{(\alpha+\varepsilon)/\beta} +\left(\frac{\delta}{\rho^{1-1/\beta}}\right)^{(\alpha-\varepsilon)n}\mathbb{E}\left(W_n^i\right)^{(\alpha-\varepsilon)/\beta} \right] \mathbb{E}\phi\left(\max_{1\leqslant r\leqslant d}\left|W_{1,0}^r-1\right|\right).
\end{eqnarray}
The assumption $\mathbb{E} \phi(Z_1^i(j))<\infty$ ($\forall1\leqslant i,j\leqslant d$) implies $\mathbb{E}\phi(\max_{1\leqslant r\leqslant d}|W_{1,0}^r-1|)<\infty$. Moreover, since $\phi(x)=x^\alpha\ell(x)$ is regularly varying with index $\alpha$, Potter's theorem implies that $$\mathbb{E}(Z_1^i(j))^{q}<\infty,\ \ 1\leqslant i,j\leqslant d, 1<q<\alpha.$$ Hence, by Theorem \ref{Thm1}, $\sup_n\mathbb{E}(W_n^i)^q<\infty.$ Since $\beta\in(1,\alpha)$, we may choose $\varepsilon>0$ sufficiently small so that $\frac{\alpha\pm\varepsilon}{\beta}\in(1,\alpha)$. Consequently, $\sup_n\mathbb{E}(W_n^i)^{(\alpha\pm\varepsilon)/\beta}<\infty$. Combining \eqref{10eq2} and \eqref{10eq3}, we obtain
\begin{align}\label{10eq4}
\mathbb{E}\phi\left(\delta^n|W_{n+1}^i-W_{n}^i|\right)\leqslant
C\left[\left(\frac{\delta}{\rho^{1-1/\beta}}\right)^{(\alpha+\varepsilon)n} +\left(\frac{\delta}{\rho^{1-1/\beta}}\right)^{(\alpha-\varepsilon)n}\right].
\end{align}
We now pass from the martingale increments to the tail $W^i-W^i_n$.
Set
\[a_{m}=(1-\delta^{-1})\delta^{-m},\quad m\geqslant 0.\]
Then $a_{m}>0$, $\sum_{m=0}^{\infty}a_{m}=1$ and $a_{m}\delta^{n+m}=(1-\delta^{-1})\delta^{n}.$
Since
$$W^{i}-W_{n}^{i}=\sum_{k=n}^{\infty}\left(W_{k+1}^{i}-W_{k}^{i}\right),$$
we have
\begin{align*}
\delta^{n}\left|W^{i}-W_{n}^{i}\right|\leqslant\sum_{m=0}^{\infty}\delta^{n}\left|W_{n+m+1}^{i}-W_{n+m}^{i}\right|
= \frac{1}{1-\delta^{-1}}\sum_{m=0}^{\infty}a_{m}\delta^{n+m}\left|W_{n+m+1}^{i}-W_{n+m}^{i}\right|.
\end{align*}
Using the convexity of $\phi$, together with $\phi(2x)\leqslant c\phi(x)$ for some $c>0$ and all $x \geqslant 0$, we have
\begin{align*}
\phi(\delta^{n}\left|W^{i}-W_{n}^{i}\right|)&\leqslant C\phi\left(\sum_{m=0}^{\infty}a_{m}\delta^{n+m}\left|W_{n+m+1}^{i}-W_{n+m}^{i}\right|\right)\\
&\leqslant C\sum_{m=0}^{\infty}a_{m}\phi\left(\delta^{n+m}\left|W_{n+m+1}^{i}-W_{n+m}^{i}\right|\right).
\end{align*}
Taking expectations and applying \eqref{10eq4}, we get
\begin{align}\label{10eq5}
\mathbb{E}\phi(\delta^{n}\left|W^{i}-W_{n}^{i}\right|)&\leqslant C\sum_{m=0}^{\infty}a_{m}\mathbb{E}\phi\left(\delta^{n+m}\left|W_{n+m+1}^{i}-W_{n+m}^{i}\right|\right)\nonumber\\
&\leqslant C\sum_{m=0}^{\infty}a_{m}\left[\left(\frac{\delta}{\rho^{1-1/\beta}}\right)^{(\alpha+\varepsilon)(n+m)} +\left(\frac{\delta}{\rho^{1-1/\beta}}\right)^{(\alpha-\varepsilon)(n+m)}\right].
\end{align}
For every $1<\delta<\rho^{1-1/\min{\{\alpha,2\}}}$, we may choose
$$\beta=\frac{1}{1-\log_\rho\left(\frac{\delta+\rho^{1-1/\min{\{\alpha,2\}}}}{2}\right)}\in(1,\alpha)$$
such that $\delta/\rho^{1-1/\beta}<1.$ Since $a_{m}=(1-\delta^{-1})\delta^{-m}\leqslant 1-\delta^{-1}$,
$$\sum_{m=0}^{\infty}a_{m}\left(\frac{\delta}{\rho^{1-1/\beta}}\right)^{(\alpha\pm\varepsilon)(n+m)}
\leqslant (1-\delta^{-1})\frac{\left(\frac{\delta}{\rho^{1-1/\beta}}\right)^{(\alpha\pm\varepsilon)n}}{1-\left(\frac{\delta}{\rho^{1-1/\beta}}\right)
^{\alpha\pm\varepsilon}}.$$
Furthermore,
\begin{align*}
\mathbb{E}\phi(\delta^{n}\left|W^{i}-W_{n}^{i}\right|)
&\leqslant C\left[\left(\frac{\delta}{\rho^{1-1/\beta}}\right)^{(\alpha+\varepsilon)n}
+\left(\frac{\delta}{\rho^{1-1/\beta}}\right)^{(\alpha-\varepsilon)n}\right]\nonumber\\
&\longrightarrow0.
\end{align*}
This proves part $(i).$

 $(ii)$ We assume that $\phi_{1/\eta}(x)$ is concave. In fact, $\phi$
is equivalent to
\[\phi_1(x)=\alpha\int_0^x t^{\alpha-1}\ell(t)\,dt,\quad x\geqslant 0.\]
Notice that
\[\phi_1'(x)=\alpha x^{\alpha-1}\ell(x) \ \text{ and } \  \ell'(x)=\ell(x)\epsilon(x)/x,\]
we have
\begin{align*}
\frac{d}{dx}\phi_1(x^{1/\eta})=\frac{\alpha}{\eta}x^{\alpha/\eta-1}\ell(x^{1/\eta}) \
\text{ and } \
\frac{d^2}{dx^2}\phi_1(x^{1/\eta})=\frac{\alpha}{\eta}x^{\alpha/\eta-2}\ell(x^{1/\eta})\left[\frac{\alpha}{\eta}-1+\frac{1}{\eta}\epsilon(x^{1/\eta})
\right].
\end{align*}
Since $\epsilon(x)\to0(x\rightarrow\infty)$, taking $a>0$ sufficiently large such that
\[\frac{\alpha}{\eta}-1+\frac{1}{\eta}\epsilon(x^{1/\eta})<0,\quad x\geqslant a,\]
we obtain that $\phi_1(x^{1/\eta})$ is concave for $x\geqslant a$. When $x<a$, set $\phi_1(x)=x^\alpha \ell(a)$
and hence, $\phi_1(x^{1/\eta})=x^{\alpha/\eta}\ell(a)$ is also concave.
Therefore, $\phi_1(x^{1/\eta})$ is concave on $[0,\infty)$. Without loss of generality, we assume that $\phi_{1/\eta}(x)$ is concave.

By an argument similar to that used in \eqref{XEQ4}, we obtain
\begin{eqnarray*}
\sum\limits_{k=1}^{\infty}\delta^{\eta n}\left|D_{n,k}^i-D_{n,k-1}^i\right|^{\eta}
\leqslant \delta^{\eta n}W_n^i \max_{1\leqslant j\leqslant d}\sum_{k=1}^{Z_n^i(j)}\frac{1}{Z_n^i(j)}\frac{|W_{k,n}^j-1|^{\eta}}{(\rho^n u(i))^{\eta-1}}.
\end{eqnarray*}
Applying Lemma \ref{lem3}, the concavity of $\phi_{1/\eta}$, and Potter's theorem,
\begin{align*}
\mathbb{E}\phi(\delta^n|W_{n+1}^i - W_n^i||\mathcal{F}_{n})&\leqslant  C\mathbb{E}\left[\phi_{1/2}\left(\sum_{k=1}^{\infty} \delta^{2n}|D_{n,k}^i - D_{n,k-1}^i|^2\right)|\mathcal{F}_{n}\right]\\
&\leqslant C\mathbb{E}\left[\phi_{1/\eta}\left(\sum_{k=1}^{\infty}\delta^{\eta n}|D_{n,k}^i - D_{n,k-1}^i|^\eta\right)|\mathcal{F}_{n}\right]\\
&\leqslant C\mathbb{E}\left[\phi_{1/\eta}\left(\delta^{\eta n}W_n^i \max_{1\leqslant j\leqslant d}\sum_{k=1}^{Z_n^i(j)}\frac{1}{Z_n^i(j)}\frac{|W_{k,n}^j-1|^{\eta}}{(\rho^n u(i))^{\eta-1}}\right)|\mathcal{F}_{n}\right]\\
&\leqslant C\sum_{j=1}^d\phi_{1/\eta}\left[\frac{\delta^{\eta n}}{\rho^{n(\eta-1)} }(W_n^i)\frac{\mathbb{E}|W_{1,0}^j-1|^{\eta}}{(u(i))^{\eta-1}}\right]\\
&\leqslant C\left[\left(\frac{\delta}{\rho^{1-1/\eta}}\right)^{(\alpha+\varepsilon)n}(W_n^i)^{\frac{\alpha+\varepsilon}{\eta}}
+\left(\frac{\delta}{\rho^{1-1/\eta}}\right)^{(\alpha-\varepsilon)n}(W_n^i)^{\frac{\alpha-\varepsilon}{\eta}}\right]
\sum_{j=1}^d\phi_{1/\eta}\left(\mathbb{E}|W_{1,0}^j-1|^{\eta}\right).
\end{align*}
Notice that $\mathbb{E}(Z_1^i(j))^{\eta}<\infty$ implies that $\mathbb{E}|W_{1,0}^j-1|^{\eta}<\infty$, and for $\epsilon>0$ small enough such that $\alpha\pm\epsilon\leqslant \eta$,
$$\mathbb{E}(W_n^i)^{\frac{\alpha\pm\varepsilon}{\eta}}\leqslant \mathbb{E}\left(1+W_n^i\mathbf{1}_{\{W_n^i>1\}}\right)
\leqslant1+\mathbb{E}(W_n^i)=2.$$
Hence, taking expectations on both sides,
\begin{align*}
\mathbb{E}\phi(\delta^n|W_{n+1}^i - W_n^i|)\leqslant  C\left[\left(\frac{\delta}{\rho^{1-1/\eta}}\right)^{(\alpha+\varepsilon)n}
+\left(\frac{\delta}{\rho^{1-1/\eta}}\right)^{(\alpha-\varepsilon)n}\right].
\end{align*}
Then following the proof of case $(i)$, we see that \eqref{10eq5} still holds true with $\eta$ replacing $\beta$, and thus
\begin{align*}
\mathbb{E}\phi(\delta^{n}\left|W^{i}-W_{n}^{i}\right|)\longrightarrow0.
\end{align*}
This proves part $(ii).$
\end{proof}
\section{Proof of Theorem \ref{Thm6-1}}
In this section, we deal with the asymptotic behavior at the critical rate $\rho^{1/2}$, and show the upper bound of the convergence rates.
\begin{proof}[Proof of Theorem \ref{Thm6-1}]
$(i)$
We first consider the case $\alpha\neq 2$. When $\alpha>2$, take $\beta=2$, whereas when $1<\alpha<2$, take $\eta=2$. In both cases, the estimate obtained in \eqref{10eq5} remains valid. Hence, at the critical value $\delta=\rho^{1/2}$,
\[
\mathbb{E}\phi\left(\delta^{n}\left|W^{i}-W_{n}^{i}\right|\right)\leqslant C\sum_{m=0}^{\infty}a_{m}.
\]
Since $\sum_{m=0}^{\infty}a_{m}=1,$
we obtain
\[
\limsup_{n\to\infty}\mathbb{E}\phi\left(\rho^{n/2}\left|W^{i}-W_{n}^{i}\right|\right)<\infty.
\]
This proves part $(i).$

$(ii)$
We next consider the case $\alpha=2$. Applying Lemma \ref{lem3}, we have
\begin{align}\label{10eq8}
\mathbb{E}\phi(\delta^n|W_{n+1}^i - W_n^i||\mathcal{F}_{n})&\leqslant  C\mathbb{E}\left[\phi_{1/2}\left(\sum_{k=1}^{\infty} \delta^{2n}|D_{n,k}^i - D_{n,k-1}^i|^2\right)|\mathcal{F}_{n}\right]\nonumber\\
&\leqslant C\mathbb{E}\left[\phi_{1/2}\left(\delta^{2 n}W_n^i \max_{1\leqslant j\leqslant d}\sum_{k=1}^{Z_n^i(j)}\frac{1}{Z_n^i(j)}\frac{|W_{k,n}^j-1|^{2}}{\rho^n u(i)}\right)|\mathcal{F}_{n}\right].
\end{align}
We first suppose that $x\rightarrow \phi(x^{1/2})$ is concave. By concavity and Potter's theorem,
\begin{align*}
\mathbb{E}\phi(\delta^n|W_{n+1}^i - W_n^i||\mathcal{F}_{n})
&\leqslant C\sum_{j=1}^d\phi_{1/2}\left[\mathbb{E}\left(\delta^{2 n}W_n^i \sum_{k=1}^{Z_n^i(j)}\frac{1}{Z_n^i(j)}\frac{|W_{k,n}^j-1|^{2}}{\rho^n u(i)}|\mathcal{F}_{n}\right)\right]\\
&= C\sum_{j=1}^d\phi_{1/2}\left[\frac{\delta^{2 n}}{\rho^n u(i)}W_n^i\mathbb{E}|W_{1,0}^j-1|^{2}\right]\\
&\leqslant C\left(\frac{\delta^{2 n}}{\rho^n u(i)}\right)^{(2+\varepsilon)/2}(W_n^i)^{(2+\varepsilon)/2}\sum_{j=1}^d\phi_{1/2}\left(\mathbb{E}|W_{1,0}^j-1|^{2}\right)\\
& \ \ \ \ \ +C\left(\frac{\delta^{2 n}}{\rho^n u(i)}\right)^{(2-\varepsilon)/2}(W_n^i)^{(2-\varepsilon)/2}\sum_{j=1}^d\phi_{1/2}\left(\mathbb{E}|W_{1,0}^j-1|^{2}\right).
\end{align*}
Since $\mathbb{E}(Z_1^i(j))^{2}<\infty$, we have $\mathbb{E}|W_{1,0}^j-1|^{2}<\infty$. Moreover, by Theorem \ref{Thm1},
$\sup_n\mathbb{E}(W_n^i)^{(2\pm\varepsilon)/2}<\infty.$ Taking expectations on both sides, we get
\begin{align*}
\mathbb{E}\phi(\delta^n|W_{n+1}^i - W_n^i|)\leqslant C\left(\frac{\delta}{\rho^{1/2}}\right)^{(2+\varepsilon)n}+C\left(\frac{\delta}{\rho^{1/2}}\right)^{(2-\varepsilon)n}.
\end{align*}
By the same argument as that used in deriving \eqref{10eq5},
\begin{align*}
\mathbb{E}\phi(\delta^{n}\left|W^{i}-W_{n}^{i}\right|)\leqslant C\sum_{m=0}^{\infty}a_{m}\left[\left(\frac{\delta}{\rho^{1/2}}\right)^{(2+\varepsilon)(n+m)} +\left(\frac{\delta}{\rho^{1/2}}\right)^{(2-\varepsilon)(n+m)}\right].
\end{align*}
Let $\delta=\rho^{1/2}$, we have \eqref{10eq6} holds.

Now suppose that $\phi_{1/2}$ is convex. Using this estimate in (\ref{10eq8}), Potter's theorem and the convexity of $\phi_{1/2},$
\begin{align*}
\mathbb{E}\phi(\delta^n|W_{n+1}^i - W_n^i||\mathcal{F}_{n})
&\leqslant C\left(\frac{\delta^{2 n}}{\rho^n u(i)}\right)^{(2+\varepsilon)/2}(W_n^i)^{(2+\varepsilon)/2}\mathbb{E}\left[\phi_{1/2}\left(\max_{1\leqslant j\leqslant d}\sum_{k=1}^{Z_n^i(j)}\frac{|W_{k,n}^j-1|^{2}}{Z_n^i(j)}\right)|\mathcal{F}_{n}\right]\\
& \ \ \ \ \ +C\left(\frac{\delta^{2 n}}{\rho^n u(i)}\right)^{(2-\varepsilon)/2}(W_n^i)^{(2-\varepsilon)/2}\mathbb{E}\left[\phi_{1/2}\left(\max_{1\leqslant j\leqslant d}\sum_{k=1}^{Z_n^i(j)}\frac{|W_{k,n}^j-1|^{2}}{Z_n^i(j)}\right)|\mathcal{F}_{n}\right]\\
&\leqslant C\left(\frac{\delta}{\rho^{1/2} }\right)^{(2+\varepsilon)n}(W_n^i)^{(2+\varepsilon)/2}\sum_{j=1}^d\mathbb{E}\phi\left(|W_{1,0}^j-1|\right)\\
& \ \ \ \ \ +C\left(\frac{\delta}{\rho^{1/2} }\right)^{(2-\varepsilon)n}(W_n^i)^{(2-\varepsilon)/2}\sum_{j=1}^d\mathbb{E}\phi\left(|W_{1,0}^j-1|\right).
\end{align*}
Taking expectations on both sides,
\begin{align*}
\mathbb{E}\phi(\delta^n|W_{n+1}^i - W_n^i|)&\leqslant C\left(\frac{\delta}{\rho^{1/2} }\right)^{(2+\varepsilon)n}\mathbb{E}(W_n^i)^{(2+\varepsilon)/2}\sum_{j=1}^d\mathbb{E}\phi\left(|W_{1,0}^j-1|\right)\\
& \ \ \ \ \ +C\left(\frac{\delta}{\rho^{1/2} }\right)^{(2-\varepsilon)n}\mathbb{E}(W_n^i)^{(2-\varepsilon)/2}\sum_{j=1}^d\mathbb{E}\phi\left(|W_{1,0}^j-1|\right)\\
&\leqslant C\left(\frac{\delta}{\rho^{1/2} }\right)^{(2+\varepsilon)n}+C\left(\frac{\delta}{\rho^{1/2} }\right)^{(2-\varepsilon)n}.
\end{align*}
Hence, let $\delta=\rho^{1/2}$, we have \eqref{10eq6} holds.

Finally, we consider the case
 $$
\ell(tx)\leqslant C\ell(x),\quad 0<t\leqslant1,\quad x>0.
$$
For $\alpha=2.$ Since $\phi(\cdot)$ is non-decreasing, applying \eqref{eq16} to the martingale difference sequence $\{\delta^n(D_{n,k}^i - D_{n,k-1}^i)\},$
\begin{align}\label{XEQ9}
\mathbb{E}\phi(\delta^n|W_{n+1}^i - W_n^i|) &\leqslant C\mathbb{E}\phi\left[\delta^n\left(\sum_{k=1}^{\infty} \mathbb{E}\left(|D_{n,k}^i - D_{n,k-1}^i|^2 \mid \mathcal{G}_{n,k-1}\right)\right)^{1/2}\right] \nonumber\\
&\ \ \ \ \ \ \ \ \ \   + C\sum_{k=1}^{\infty} \mathbb{E}\phi(\delta^n|D_{n,k}^i - D_{n,k-1}^i|).
\end{align}
We first estimate the first term on the right hand side. For $r\neq s,$
\begin{align*}
\mathbb{E}\left[(W_{k,n}^r - 1)(W_{k,n}^s - 1) \mid \mathcal{G}_{n,k-1}\right] = \mathbb{E}(W_{k,n}^r - 1)\mathbb{E}(W_{k,n}^s - 1) = 0.
\end{align*}
Hence,
\begin{align*}
&\mathbb{E}\left(|D_{n,k}^i - D_{n,k-1}^i|^2 \mid \mathcal{G}_{n,k-1}\right)\\
&=\mathbb{E}\left[\left|\sum_{r=1}^{d} \frac{u(r)}{\rho^n u(i)}(W_{k,n}^r - 1)\mathbf{1}_{\{Z_n^i(r)\geqslant k\}}\right|^2 \Bigg| \mathcal{G}_{n,k-1}\right]\\
&=\sum_{r=1}^{d} \frac{u(r)^2}{\rho^{2n}u(i)^2}\mathbf{1}_{\{Z_n^i(r)\geqslant k\}} \mathbb{E}|W_{k,n}^r - 1|^2\\
& \ \ \ \ \ \ \ \ +2\sum_{1\leqslant r<s\leqslant d} \frac{u(r)u(s)}{\rho^{2n}u(i)^2}\mathbf{1}_{\{Z_n^i(r)\geqslant k\}}\mathbf{1}_{\{Z_n^i(s)\geqslant k\}}\mathbb{E}\left[(W_{k,n}^r - 1)(W_{k,n}^s - 1) \mid \mathcal{G}_{n,k-1}\right]\\
&=\sum_{r=1}^{d} \frac{u(r)^2}{\rho^{2n}u(i)^2}\mathbf{1}_{\{Z_n^i(r)\geqslant k\}} \mathbb{E}|W_{1,0}^r - 1|^2.
\end{align*}
Summing over $k,$
\begin{align*}
\sum_{k=1}^{\infty} \mathbb{E}\left(|D_{n,k}^i - D_{n,k-1}^i|^2 \mid \mathcal{G}_{n,k-1}\right)
&\leqslant\sum_{k=1}^{\infty}\sum_{r=1}^{d} \frac{u(r)^2}{\rho^{2n}u(i)^2}\mathbf{1}_{\{Z_n^i(r)\geqslant k\}} \mathbb{E}|W_{1,0}^r - 1|^2\\
&=\sum_{r=1}^{d} \frac{u(r)^2Z_n^i(r)}{\rho^{2n}u(i)^2} \mathbb{E}|W_{1,0}^r - 1|^2\\
&\leqslant C\rho^{-n}W_n^i.
\end{align*}
By the condition  $\ell(tx)\leqslant C\ell(x), 0<t\leqslant1, x>0,$ we have
\begin{align*}
\phi[(W_n^i)^{1/2}]&=\phi[(W_n^i)^{1/2}]\mathbf{1}_{\{W_n^i\geqslant1\}}
+\phi[(W_n^i)^{1/2}]\mathbf{1}_{\{W_n^i<1\}}\\
&\leqslant W_n^i\ell[(W_n^i)^{1/2}]\mathbf{1}_{\{W_n^i\geqslant1\}}+\phi(1)\\
&=W_n^i\ell[(W_n^i)^{-1/2}(W_n^i)]\mathbf{1}_{\{W_n^i\geqslant1\}}+\phi(1)\\
&\leqslant CW_n^i\ell(W_n^i)\mathbf{1}_{\{W_n^i\geqslant1\}}+\phi(1)\\
&\leqslant C\phi(W_n^i)+\phi(1).
\end{align*}
By Theorem \ref{Thm2}, $\mathbb{E}\phi\left(Z^i_1(j)\right)<\infty$ implies that $\mathbb{E}\phi(\widehat{W}^i)<\infty$. Hence, at the critical value $\delta=\rho^{1/2},$
\begin{align}\label{10eq9}
\mathbb{E}\phi\left[\delta^n\left(\sum_{k=1}^{\infty}\mathbb{E}\left(|D_{n,k}^i - D_{n,k-1}^i|^2 \mid \mathcal{G}_{n,k-1}\right)\right)^{1/2}\right]
\leqslant\mathbb{E}\phi\left[C(W_n^i)^{1/2}\right]
\leqslant C\left[\mathbb{E}\phi(W_n^i)+\phi(1)\right]
<\infty.
\end{align}
Next, we estimate the second term on the right hand side of \eqref{XEQ9}. Since $\phi$ is convex on $[0,\infty)$ with $\phi(2x) \leqslant c\phi(x), x>0$, there exists a constant $C$ such that
$$\phi(x_1 + \cdots + x_d) \leqslant C\sum_{r=1}^{d} \phi(x_r), \ \ x_r\geqslant0.$$
By the triangle inequality and the monotonicity of $\phi$, we have
\begin{align*}
\phi\left(\delta^n\left|D_{n,k}^i-D_{n,k-1}^i\right|\right)
&\leqslant\phi\left(\sum_{r=1}^d \frac{\delta^nu(r)}{\rho^n u(i)}\left|W_{k,n}^r-1\right| \mathbf{1}_{\left\{Z_n^i(r) \geqslant k\right\}}\right)\\
&\leqslant C\sum_{r=1}^d\phi\left(\frac{\delta^nu(r)}{\rho^n u(i)}\left|W_{k,n}^r-1\right| \mathbf{1}_{\left\{Z_n^i(r) \geqslant k\right\}}\right)\\
&\leqslant C\sum_{r=1}^d\phi\left(C\frac{\delta^n}{\rho^n}\left|W_{k,n}^r-1\right| \right)\mathbf{1}_{\left\{Z_n^i(r) \geqslant k\right\}}.
\end{align*}
Hence,
\begin{align*}
\sum_{k=1}^{\infty} \mathbb{E}\phi\left(\delta^n|D_{n,k}^i - D_{n,k-1}^i|\right)
&\leqslant C\sum_{k=1}^{\infty}\mathbb{E}\sum_{r=1}^d\phi\left(C\frac{\delta^n}{\rho^n}\left|W_{k,n}^r-1\right| \right)\mathbf{1}_{\left\{Z_n^i(r) \geqslant k\right\}}\\
&\leqslant C\sum_{r=1}^{d} \mathbb{E}\sum_{k=1}^{Z_n^i(r)} \phi\left(C\frac{\delta^n}{\rho^n}|W_{k,n}^r - 1|\right)\\
&\leqslant C\sum_{r=1}^{d}\rho^n\mathbb{E}\phi\left(\frac{\delta^n}{\rho^n}|W_{1,0}^r - 1|\right).
\end{align*}
By the condition  $\ell(tx)\leqslant C\ell(x), 0<t\leqslant1, x>0,$ we have
\[
\phi(tx)=t^{2}x^{2}\ell(tx)\leqslant C t^{2}\phi(x),\qquad 0<t\leqslant 1.
\]
Thus, at $\delta=\rho^{1/2},$
\begin{align*}
\sum_{k=1}^{\infty} \mathbb{E}\phi\left(\delta^n|D_{n,k}^i - D_{n,k-1}^i|\right)
\leqslant C\sum_{r=1}^{d}\mathbb{E}\phi\left(|W_{1,0}^r - 1|\right)<\infty.
\end{align*}
Combining (\ref{XEQ9}) and (\ref{10eq9}), we obtain
\begin{align*}
\sup_n\mathbb{E}\phi(\rho^{n/2}|W_{n+1}^i - W_n^i|)<\infty.
\end{align*}
By the same argument as that used in deriving \eqref{10eq5}, let $\delta=\rho^{1/2}$, we have \eqref{10eq6} holds.
\end{proof}

\section{Proof of Theorem \ref{Thm7}}
\begin{proof}[Proof of Theorem \ref{Thm7}]
By $\eqref{eq15}$, for all $n\geqslant0$ and $1\leqslant i\leqslant d,$
\begin{align}\label{Xiaoeq19}
W^i-W_n^i=\sum_{j=1}^d \frac{u(j)}{\rho^{n}u(i)} \sum_{l=1}^{Z_n^i(j)}\left(W_{l, n, \infty}^j-1\right),
\end{align}
where the random variables $W_{l, n, \infty}^j$, indexed by $l\geqslant1,n \geqslant 0$ and $1 \leqslant j \leqslant d$, are independent of each other and independent of $Z_n^i$, and have the same distribution as $W^j.$ Under the assumed moment conditions, we have $\mathbb{E}(W_{l, n, \infty}^j-1)=\mathbb{E}(W^j-1)=0.$
For the fixed type $1\leqslant j\leqslant d,$ we decompose \eqref{Xiaoeq19} as
\begin{align*}\label{XEQ13}
W^i-W_n^i=\frac{u(j)}{\rho^nu(i)}\sum_{l=1}^{Z_n^i(j)}\left(W_{l, n,\infty}^j-1\right)+\sum_{\substack{1 \leqslant s \leqslant d \\
s \neq j}} \frac{u(s)}{\rho^nu(i)} \sum_{l=1}^{Z_n^i(s)}\left(W_{l, n,\infty}^s-1\right).
\end{align*}
Let
\begin{align*}
\mathcal{H}_{n,j}=\sigma(\mathcal{F}_n,W_{l,n,\infty}^j:l\geqslant1).
\end{align*}
By the independence between different types, for each $s \neq j,$
\begin{equation*}\label{XEQ14}
\mathbb{E}\left[\sum_{l=1}^{Z_n^i(s)}\left(W_{l, n,\infty}^s-1\right) \mid \mathcal{H}_{n,j}\right]
=\sum_{l=1}^{Z_n^i(s)}\mathbb{E}\left[\left(W_{l, n,\infty}^s-1\right) \mid \mathcal{F}_n\right]
=\sum_{l=1}^{Z_n^i(s)}\mathbb{E}\left[\left(W^s-1\right)\right]=0.
\end{equation*}
Consequently,
$$\mathbb{E}(W^i - W_n^i \mid \mathcal{H}_{n,j}) = \frac{u(j)}{\rho^n u(i)} \sum_{l=1}^{Z_n^i(j)} (W_{l,n,\infty}^j - 1).$$
Since $\phi$ is increasing and convex, the function $x \longmapsto \phi\left(\delta^n|x|\right)$ is convex. By Jensen's inequality,
\begin{align*}
\mathbb{E}\left[\phi\left(\delta^{n}\left|W^{i}-W_{n}^{i}\right|\right) \mid \mathcal{H}_{n, j}\right]
&\geqslant \phi\left(\delta^{n}\left|\mathbb{E}\left(W^{i}-W_{n}^{i} \mid \mathcal{H}_{n, j}\right)\right|\right).\\
&= \phi\left(\frac{\delta^{n} u(j)}{\rho^{n} u(i)}\left|\sum_{l=1}^{Z_{n}^{i}(j)}\left(W_{l, n, \infty}^{j}-1\right)\right|\right).
\end{align*}
Taking expectations yields
\begin{align}\label{XEQ20}
\mathbb{E}\phi\left(\delta^{n}\left|W^{i}-W_{n}^{i}\right|\right)\geqslant\mathbb{E}\phi\left(\frac{\delta^{n} u(j)}{\rho^{n} u(i)}\left|\sum_{l=1}^{Z_{n}^{i}(j)}\left(W_{l, n, \infty}^{j}-1\right)\right|\right).
\end{align}
For fixed $n$ and $j$, the sequence $\{\sum_{l=1}^k(W_{l, n,\infty}^j-1) \mathbf{1}_{\left\{Z_n^i(j) \geqslant l\right\}}\}_{k \geqslant 1}$ forms a martingale with respect to the filtration $\sigma(\mathcal{F}_n,W_{l,n,\infty}^j:1\leqslant l\leqslant k)$. Applying Lemma \ref{lem3}(ii) gives
\begin{align}\label{XEQ21}
\mathbb{E}\phi\left(\frac{\delta^{n}u(j)}{\rho^{n}u(i)}\left|\sum_{l=1}^{Z_{n}^{i}(j)}\left(W_{l,n,\infty}^{j}-1\right)\right|\right)
\geqslant C\mathbb{E}\phi\left(\frac{\delta^{n}u(j)}{\rho^{n}u(i)}\left[\sum_{l=1}^{Z_{n}^{i}(j)}\left|W_{l,n,\infty}^{j}-1\right|^{2}\right]^{1/2}\right).
\end{align}
Combining \eqref{XEQ20} and \eqref{XEQ21}, and then using Jensen's inequality, we obtain
\begin{align*}
\mathbb{E}\phi\left(\delta^{n}\left|W^{i}-W_{n}^{i}\right|\right)
&\geqslant C\mathbb{E}\phi\left(\frac{\delta^{n}u(j)}{\rho^{n}u(i)}\left[\sum_{l=1}^{Z_{n}^{i}(j)}\left|W_{l,n,\infty}^{j}-1\right|^{2}\right]^{1/2}\right)\\
&\geqslant
C\phi\left(\frac{\delta^{n}u(j)}{\rho^{n}u(i)}\mathbb{E}\left[\sum_{l=1}^{Z_{n}^{i}(j)}\left|W_{l,n,\infty}^{j}-1\right|^{2}\right]^{1/2}\right).
\end{align*}
By inequality $(x_1+\cdots+x_n)^{1/2} \geqslant n^{-1/2}(x_1^{1/2}+\cdots+x_n^{1/2})$ for $x_1, \cdots, x_n \geqslant 0$, we deduce
\begin{align*}
\mathbb{E}\left(\sum_{l=1}^{Z_n^i(j)}\left|W_{l, n,\infty}^j-1\right|^2\right)^{1 / 2}
&=\mathbb{E}\left[\mathbb{E}\left(\sum_{l=1}^{Z_n^i(j)}\left|W_{l, n,\infty}^j-1\right|^2\right)^{1 / 2}|\mathcal{F}_n\right]\nonumber\\
&\geqslant\mathbb{E}\left[Z_n^i(j)^{-1/2}\mathbb{E}\left(\sum_{l=1}^{Z_n^i(j)}\left|W_{l, n,\infty}^j-1\right||\mathcal{F}_n\right)\right]\nonumber\\
&=\mathbb{E}Z_n^i(j)^{1/2}\mathbb{E}|W^j-1|.
\end{align*}
It follows that
\begin{equation}\label{XEQ15}
\mathbb{E}\phi\left(\delta^{n}\left|W^{i}-W_{n}^{i}\right|\right)
\geqslant C\phi\left(\frac{\delta^{n}u(j)}{\rho^{n}u(i)}\mathbb{E}Z_{n}^{i}(j)^{1/2}\mathbb{E}\left|W^{j}-1\right|\right).
\end{equation}
By the Kesten-Stigum theorem,
\[
\frac{Z_{n}^{i}(j)}{\rho^{n}v(j)} \to u(i)W^{i} \ \  \text{a.s.}
\]
and hence
\[
\rho^{-n/2}Z_{n}^{i}(j)^{1/2} \to \bigl(v(j)u(i)W^{i}\bigr)^{1/2} \ \ \text{a.s.}
\]
Fatou's lemma therefore gives
\begin{equation}\label{XEQ16}
\begin{aligned}
\liminf_{n\to\infty}\rho^{-n/2}\mathbb{E}Z_{n}^{i}(j)^{1/2}\geqslant \bigl(v(j)u(i)\bigr)^{1/2}\mathbb{E}\bigl(W^{i}\bigr)^{1/2}> 0.
\end{aligned}
\end{equation}
Moreover, by the non-triviality assumption, we have $\mathbb{E}|W^j-1|>0$. Consequently, by \eqref{XEQ15} and \eqref{XEQ16}, for all sufficiently large $n,$
\begin{align}\label{liang-1}
\mathbb{E}\phi\left(\delta^{n}\left|W^{i}-W_{n}^{i}\right|\right)
&\geqslant C\phi\left(C\frac{u(j)}{u(i)}\left(\frac{\delta}{\rho^{1/2}}\right)^n\right)\geqslant C\phi\left(C\left(\frac{\delta}{\rho^{1/2}}\right)^n\right).
\end{align}

When $\delta>\rho^{1/2}$, since $\phi(x)\rightarrow\infty$ as $x\rightarrow\infty$, \eqref{liang-1} and the fact that $(\delta/\rho^{1/2})^n\rightarrow\infty$ yield $$\mathbb{E}\phi\left(\delta^{n}\left|W^{i}-W_{n}^{i}\right|\right)\rightarrow\infty.$$
When $\delta=\rho^{1/2}$, \eqref{liang-1} gives
$$\liminf_{n\rightarrow\infty}\mathbb{E}\phi\left(\delta^{n}\left|W^{i}-W_{n}^{i}\right|\right)\geqslant C\phi(C)>0.$$
This completes the proof.
\end{proof}

\noindent\textbf{Acknowledgements} The work has been partially supported by the National Natural Science Foundation of China (Grants nos.\ 11971063, 11601019, 12271062), and by the Natural Science Foundation of Beijing, China (Grant No. 1252005).

\end{document}